\documentstyle[12pt]{article}
\newcommand{\be}{\begin{equation}}
\newcommand{\ee}{\end{equation}}
\newcommand{\og}{({\cal K}, ~{\cal L}, ~{\cal H}, ~{\cal J})}
\newcommand{\al}{\alpha}
\newcommand{\pp}{{\cal PP}'}
\newcommand{\bet}{\beta}
\newcommand{\g}{\gamma}
\newcommand{\om}{\omega}

\newcommand{\G}{\Gamma}

\newcommand{\La}{\Lambda}

\newcommand{\crs}{{\cal Q}}
\newcommand{\M}{{\cal M}}
\newcommand{\Pen}{{\cal P}}

\newcommand{\la}{\lambda}
\newcommand{\eps}{\varepsilon}

\newcommand{\dz}{\wedge}
\newcommand{\C}{{\bf C }}

\newcommand{\ba}{\begin{array}}
\newcommand{\ea}{\end{array}}
\newcommand{\beq}{\begin{eqnarray}}
\newcommand{\eeq}{\end{eqnarray}}
\newtheorem{lm}{Lemma}[section]
\newtheorem{sss}{Theorem}[section]
\newtheorem{pr}{Proposition}[section]
\newtheorem{co}{Corollary}[section]
\newtheorem{exa}{Example}[section]
\newtheorem{deff}{Definition}[section]
\newcommand{\bd}{\begin{deff}}
\newcommand{\ed}{\end{deff}}
\newcommand{\bl}{\begin{lm}}
\newcommand{\el}{\end{lm}}
\newcommand{\bp}{\begin{pr}}
\newcommand{\ep}{\end{pr}}
\newcommand{\bt}{\begin{sss}}
\newcommand{\et}{\end{sss}}
\newcommand{\bc}{\begin{co}}
\newcommand{\ec}{\end{co}}
\newcommand{\bex}{\begin{exa}}
\newcommand{\eex}{\end{exa}}

\newcommand{\der}{\mbox{d}}
\begin{document}

\thispagestyle{empty}

\title {TWISTOR BUNDLES, EINSTEIN EQUATIONS AND REAL STRUCTURES
\footnote{Research supported 
by SISSA, SFB 288 and 
Polish Committee for Scientific Research (KBN) Grant nr 2 P302 112 7.}} 
     \author{ 
{\sc Pawe\l  ~Nurowski}\\
{\small {\it Department of Mathematical Methods in Physics}}\\
{\small{\it Faculty of Physics, Warsaw University}}\\
{\small {\it ul. Ho\.za 74, Warszawa, Poland}}\\ 
{\small {\it e-mail nurowski@fuw.edu.pl}}\\ 
            }
\maketitle
\begin{abstract}
We consider ${\bf S}^2$ bundles $\Pen$ and $\Pen'$ of totally null planes of 
maximal dimension over a 4-dimensional manifold equipped with a Weyl or 
Riemannian geometry. The fibre product $\pp$ of $\Pen$ and $\Pen'$ is found 
to be appropriate for the encoding of both the selfdual and the Einstein-Weyl 
equations 
for the 4-metric. This encoding 
is realized  in terms of the properties of certain well defined geometrical 
objects on $\pp$. The formulation is suitable for complex-valued metrics 
and unifies results for all three 
possible real signatures. In the purely Riemannian positive definite case 
it implies the existence of a natural almost hermitian structure on $\pp$ 
whose integrability conditions correspond to the self-dual Einstein 
equations of the 4-metric. All Einstein equations for the 4-metric are 
also encoded in the properties of this almost hermitian structure on $\pp$. 

\end{abstract}
\newpage

\section{Introduction}
The natural appearance of complex coordinates in the Robinson-Trautman 
\cite{bi:RTm} class of metrics was one of the first signs that 
complex geometric methods may be important in General Relativity. 
The proper understanding of this fact led to the introduction of 
CR-geometric concepts to the Einstein 
theory \cite{bi:LNT,bi:Pen3,bi:RT1,bi:RT,bi:Taf}. Penrose's twistor programme 
was also partially motivated by this result.\\
In this paper we are concerned with a twistor theory over a 4-dimensional 
manifold. Given a metric on such a manifold the problem of  
encoding the Einstein equations on the associated twistor bundle arises. 
Since Penrose's original paper \cite{bi:Pen2} several 
attempts to construct an encoding has been undertaken 
\cite{bi:New,bi:Sal,bi:Spar}. 
In particular, in the case of a positive definite metric, 
S. Salamon in Ref. \cite{bi:Sal} 
used well defined differential forms on the twistor bundle and showed that 
the vanishing of certain differentials corresponded to the anti-self-dual 
Einstein equations on the base manifold. Our approach in Ref. 
\cite{bi:optical} is very 
much in the spirit of Salamon. There we studied  
natural forms on the twistor bundle 
in the Lorentzian case. Our analysis was more complex than in the 
positive definite case since we had to deal with directions of forms 
rather than forms themselves. We showed that if our forms 
satisfied certain well defined differential conditions on the twistor bundle 
then the Ricci tensor of the base metric was traceless.  
Thus, in that case, we succeeded in encoding the full set of Einstein 
equations, without restricting to 
anti-selfdual metrics. In this paper we extend 
results of \cite{bi:optical} to 4-manifolds with complex 
valued metrics or real metrics of signature (+,+,+,+) or (+,+,--,--). 
The proposed approach unifies all the signatures and applies also to Weyl 
geometries.\\
Section 2 presents a short reformulation of the main results of 
\cite{bi:optical}. It will be useful for generalizations of the results 
to other signatures. \\
Section 3 describes an analogy between hermitian and 
optical geometries (see Theorem 3.1). We show there that 
notions such as that of a null direction in 4-dimensional Lorentzian geometry 
and 
an almost hermitian structure in the case of a positive definite metric have 
a unified description in terms of totally null planes of maximal dimension in 
the 
complexification of the tangent space. It turns out that the integrability 
conditions for both almost hermitian geometries and optical geometries 
associated with null rays have a unified description in terms of associated 
fields of maximal totally null planes. These conditions are given by equation 
(\ref{eq:integr}).\\
Section 4 gives necessary information about Weyl geometries. We recall that 
such geometries 
are given in terms of a class of pairs $(g, A)$ where $g$ is a metric 
and $A$ is a 1-form on a manifold. Two pairs $(g, A)$ and 
$(g',A')$ are in the same class iff $g'={\rm e}^{2\varphi}g$ 
and $A'=A-2\der\varphi$.\\ 
In Section 5 we study 4-manifolds equipped with 
Weyl geometries. We consider Weyl geometries in which both $g$ and $A$ may 
be complex-valued. However, we do not exclude situations in which $g$ and 
$A$ are real. Purely metric situations $A=0, \varphi=0$ are also 
not excluded in our analysis. Given a 4-manifold $\cal M$ and a Weyl geometry 
$(g,A)$ on it we consider a space $\Pen$ of all selfdual totally null 2-planes 
in the complexification of T$\cal M$. This space is an {\bf S}$^2$ bundle over 
$\cal M$. An analogous bundle $\Pen'$ of all anti-selfdual maximal totally 
null spaces is also considered there. Given $\Pen$ and $\Pen'$ we also study  
their fibre product space $\pp$. This is a bundle over $\cal M$ with  
typical fiber ${\bf S}^2\times{\bf S}^2$. We call $\Pen$, $\Pen'$ and $\pp$ 
twistor bundles (Section 5.1).\\
The rest of Section 5 is devoted to studies of natural geometric structures 
that exist on twistor bundles. 
In particular, we find that $\cal P$ has the following geometrical features. 
T$\Pen$ splits naturally into a vertical and horizontal part. One can also 
naturally define on $\Pen$ a spin 
connection 1-form, a class of metrics $\tilde{g}$,  
a canonical field of horizontal 2-planes, and two distinguished fields of 
3-planes 
which are totally null in any metric from the class $\tilde{g}$. We 
also find a way of writing certain differential equations on $\cal P$ 
that have geometrical meaning. Analogous structures are also defined on 
$\Pen'$ 
(Section 5.2). Another set of geometrical objects is naturally defined on 
$\pp$ (Section 5). There we find a natural split of the tangent bundle 
into vertical and horizontal parts. This enables a canonical field of 
horizontal 3-planes to be defined on $\pp$. There is a nice geometry 
associated 
with these which, in particular, distinguishes a certain field of 
(in general complex) 
directions. This field is null in the naturally defined class of metrics 
on $\pp$. It is used to define a canonical 1-form and eight 
distinguished fields of 4-planes that are totally null in all the 
natural metrics on $\pp$.\\
Section 6 deals with the integrability conditions of the structures 
defined on twistor bundles. Using them we find a way of encoding 
(anti)selfduality equations for the 4-metric on $\Pen$ and $\Pen'$ 
(Theorems 6.4, 6.5). This gives a Weyl-geometric generalizations of 
the Atiyah-Hitchin-Singer \cite{bi:AHS} and the Penrose \cite{bi:Pen3} 
theorems. Use of the natural structures on $\pp$ enables the 
(anti)selfdual Einstein-Weyl equations to be encoded there. This encoding is 
done by means of the integrability conditions of one of the eight naturally 
defined fields of 
maximal totally null planes on $\pp$ (Corollary 6.1, Theorem 6.3). 
Another results of this Section are included in Theorems 6.1 and 6.2. They 
provide a description of the Einstein-Weyl 
equations (without restriction to selfdual metrics) on $\pp$. 
In the purely metric case $A=0,\varphi=0$ they give a signature independent 
formulation of the Einstein equations $R_{\mu\nu}=\la g_{\mu\nu}$ on $\pp$.\\
Section 7 is concerned with the interpretations of the results of Section 6 in 
the case of real Weyl geometries. If the 4-metric has positive definite or 
neutral signature 
(Section 7.2), then the main results are included in Theorem 7.1. This, in 
particular, states that there is a preferred almost hermitian 
structure on $\pp$ 
whose integrability conditions are equivalent to the selfdual 
Einstein-Weyl equations for the Weyl geometry. This result even in the purely 
metric 
case goes a bit beyond the Atiyah-Hitchin-Singer theorem. We are able 
to encode both the selfduality and Einstein equations in the integrability 
conditions of natural almost hermitian structure on $\pp$. 
Section 7.2 also includes geometrical interpretation of 
the full set of Einstein-Weyl 
(Einstein, in pure metric case) equations on $\pp$. This is given 
by Theorem 7.2. It uses one of the eight distinguished almost 
hermitian structures $\cal J$ on $\pp$ to analyze the decomposition of the 
differential of the naturally defined spin connection 1-form on $\pp$. It 
turns out that the Einstein-Weyl equations for the Weyl geometry are 
equivalent to the fact that this differential has zero intersection with 
the T$^{*(0,2)}$ space of 2-forms, where T$^{*(0,2)}$ is defined with respect 
to $\cal J$. Section 7.3 deals  
with the Lorentzian case of the 4-metric. The main result is summarized 
in Theorem 7.3. It, in particular, states that in the purely metric case one 
can associate a natural 7-dimensional CR-structure with the Minkowski 
4-metric. 
The end of Section 7 explains why in Lorentzian case it 
suffices to work on $\Pen$ or $\Pen'$. The results of Section 2 are then 
reobtained in terms of bundles of maximal totally null planes 
(Theorems 7.4 and 7.5) .\\ 
\section{Summary of the Lorentzian case}
To make the present paper self 
contained we briefly recall our 
twistor formulation of the Einstein equations in Lorentzian case 
\cite{bi:optical}.

Let $\cal M$ be a 4-dimensional oriented manifold  equipped with a
Lorentzian metric $g$. It is convenient to introduce a null tetrad 
$(m,~\bar{m},~l,~k)$ on $\M$ with a cotetrad 
$(\theta^i)=(\theta^1,~\theta^2,~\theta^3,~\theta^4)=(M,~\bar{M},~L,~K)$ 
so that 
\be
g=g_{ij}\theta^i\theta^j=M\bar{M}-LK,
\label{eq:met}
\ee
where $\theta^i\theta^j=\frac{1}{2}
(\theta^i\otimes \theta^j+\theta^j\otimes \theta^i)$. 
Consider set {\bf S}$_x$ of all  
null directions outgoing from a given point $x\in\cal M$. This set is 
topologically a
sphere (the celestial sphere of an observer situated at $x$). The points 
of this sphere can be parameterized by a complex number $z$ belonging 
to the Argand plane ${\bf C}\cup\{\infty\}$. 
A direction associated with  $z\neq\infty$ is generated by a vector
\be
k(z)= k +z\bar{z}l-z m -\bar{z}\bar{m}\label{eq:kz}.
\ee 
With $z=\infty$ we associate a direction generated 
by vector $l$. Conversely, any null direction from $x$ is either 
parallel to the vector $l$ or can be represented by the unique null 
vector $k(z)$ such that $g(k(z),l)=-1$. It follows that $k(z)$ has 
necessarily the form (\ref{eq:kz}), and that it defines a certain 
$z\in\bf C$. If 
a direction is parallel to $l$ we associate with it $z=\infty$.\\
We define a fiber bundle  
${\cal P}=\displaystyle{\bigcup_{x\in\cal M}}{\bf S}_x$ 
over $\cal M$, so  that its fibers are 2-dimensional spheres ${\bf S}_x$. 
The anonical projection $\pi:{\cal P}\rightarrow{\cal M}$ is defined by  
$\pi({\bf S}_x)=x$. We will call the bundle $\cal P$ 
``Penrose's twistor space'',  or the `twistor bundle'. This bundle possesses
quite a reach family of well defined geometrical objects, which  
collectively form the 
so called `optical geometry' \cite{bi:ATOG,bi:T2}. Here we recall only 
those objects that are relevant in the present paper (see \cite{bi:optical} 
for details).
\begin{itemize}
\item[{\it i).}] Levi-Civita connection associated with the metric $g$ on $\M$ 
distinguishes a horizontal space in T$\cal P$. A vertical space  
consists, by definition, of vectors tangent to the fibers. In this way, at any 
point $p\in \cal P$ 
we have a natural splitting of its tangent space onto a direct sum 
T$_p{\cal P}=V_p\oplus H_p$, where $H_p$ is 4-dimensional horizontal space 
and $V_p$ is a 2-dimensional vertical space. The vertical space 
$V_p$ is identical 
with a tangent space to a certain point on the 2-dimensional sphere. Thus 
$V_p$ 
has a natural complex structure related to the complex structure on 
${\bf S}^2$. Complexification of $V_p$ splits it into eigenspaces $V^+_p$ and 
$V^-_p$ with respect to this complex structure. We have a horizontal 
lift $\tilde{v}$ of any 
vector $v$ from $\pi (p)\in\M$ to $\Pen$. This is a vector 
$\tilde{v}$ such that at $p$ $\tilde{v}\in H_p$ and $\pi_* (\tilde{v})=v$.
\item[{\it ii).}] A Lorentzian metric $\tilde{g}$ can be defined on $\cal P$ 
by the requirements that
\begin{itemize}
\item[{\it a)}] a scalar product of any two horizontal vectors is defined in
$\tilde{g}$ in terms of the scalar product in $g$ of their push 
forwards to 
$\cal M$,
\item[{\it b)}] a scalar product of any two vertical vectors in $\tilde{g}$ is 
equal to their scalar product in the natural metric on 2-dimensional 
sphere (this is consistent since vertical vectors can be considered 
tangent vectors to ${\bf S}^2$),
\item[{\it c)}] any two vectors such that one is horizontal and the other 
is vertical are orthogonal in $\tilde{g}$.
\end{itemize}
\item[{\it iii).}] There is a natural congruence of lines on $\Pen$ which is
tangent to the horizontal lifts 
of null directions from $\cal M$. It is defined by the following 
receipee. Take a null vector $k$ at $x\in\cal M$. This represents 
a certain null direction $p(k)$ outgoing from $x$. Correspondingly, 
this defines a point $p=p(k)$ in the fiber $\pi^{-1}(x)$. Lift $k$ 
horizontally to $p$. This defines $\tilde{k}$ which generates a 
certain direction 
outgoing from $p\in\cal P$. Repeating this procedure for all 
directions outgoing from $x\in\cal M$ we attacha unique direction to any point 
of $\pi^{-1}(x)$. If we do it for all points of 
$\cal M$, we define a field of directions on 
$\cal P$ which, according to its construction and properties of 
$\tilde{g}$, is null. Integral curves of this field form the desired null 
congruence. This congruence is called the null spray on $\cal P$ 
\cite{bi:Spar}.
\end{itemize}

Let $X$ 
be any non vanishing vector field tangent to the null spray on 
$\cal P$. Let $\Lambda_L$ be a real 1-form on $\cal P$ defined by 
$\Lambda_L=\tilde{g}(X)$. Since $X$ is defined up to a multiplication of 
a non vanishing real function on $\Pen$ then also $\La_L$ is specified up 
to a multiplication by a 
real non vanishing function $u$ on $\cal P$ 
\be
\Lambda_L\rightarrow\Lambda'_L=u\Lambda_L. \label{eq:trla}
\ee
One associates another 1-form with the horizontal space in $\cal P$. This is 
such a complex 1-form $E_L$ 
on $\cal P$ that satisfies  
$i)$ $E_L(H_p)=E_L(V^-_p)=0$ and $ii)$
$E_L\wedge\bar{E}_L\neq 0$ at any point $p\in\cal P$. $E_L$ is also 
defined up to a multiplication 
by a non vanishing complex function $h$ on $\cal P$ 
\be
E_L\rightarrow E'_L=h E_L.\label{eq:tre}
\ee 
It is easy to see that 
the metric $\tilde{g}$ on $\cal P$ can be expressed as
$$\tilde{g}=2(h\bar{h}E_L\bar{E}_L+\Lambda_L T + F_L\bar{F}_L)$$
with some 1-forms $T$ (real) 
and $F_L$ (complex) on $\cal P$. The above expression can be considered 
a definition of the form $F_L$. It is given up to transformations 
\be
F_L\rightarrow F'_L=\mbox{e}^{i\phi}F_L+p\Lambda_L,\label{eq:trf}
\ee
where $\phi$ (real) and $p$ (complex) are some functions on $\Pen$.\\
It follows that in the ordered null cotetrad 
$(\theta^i)$ of (\ref{eq:met}) the forms $\Lambda_L$, $F_L$ and $E_L$ can be 
represented by 
\be
\Lambda_L=-L-z\bar{z}K-z\bar{M}-\bar{z}M,\label{eq:lal}
\ee
\be
F_L=M+zK,\label{eq:fl}
\ee
\be
E_L=dz+\gamma^3_{~2}+z(\gamma^1_{~1}-\gamma^4_{~4})+z^2\gamma^2_{~3},
\label{eq:el}
\ee
where $z$ is the same as in (\ref{eq:kz}) and $\gamma^i_{~j}$ are 
Levi-Civita connections 1-forms associated with the metric $g$ in the cotetrad 
$(\theta^i)$.

Although the above forms are only defined up to transformations 
(\ref{eq:trla})-(\ref{eq:trf}) one can use them to write down some well 
defined equations on $\Pen$. The following  
equations, invariant under 
transformations (\ref{eq:trla})-(\ref{eq:trf}), are of particular interest.\\
\beq
&\der\Lambda_L\dz\Lambda_L\dz F_L\dz E_L=0,\label{eq:wc}\\
&\der F_L\dz\Lambda_L\dz F_L\dz E_L=0,\label{eq:wc'}\\
&\der E_L\dz\Lambda_L\dz F_L\dz E_L=0,\label{eq:weyll}\\
&\der E_L\dz\Lambda_L\dz \bar{F}_L\dz E_L=0.\label{eq:ricl}
\eeq 
Note that in equation (\ref{eq:ricl}) a form $\bar{F}_L$, which is a 
complex conjugate of $F_L$, appears.

\noindent
Since any of the above equations is invariant under 
(\ref{eq:trla})-(\ref{eq:trf}) we can use $\Lambda_L$, $F_L$ and $E_L$ in a 
particular representation  (\ref{eq:lal})-(\ref{eq:el}) 
to analyze them.
It is a matter of a straightforward but lengthy calculation to 
arrive at the following 
theorem.
\bt~\\
\vspace{-.8truecm}
\begin{itemize}
\item[1).] Equations (\ref{eq:wc})-(\ref{eq:wc'}) are identically satisfied on
$\Pen$.
\item[2).] Equation (\ref{eq:weyll}) is satisfied everywhere on $\Pen$ if and 
only if the metric (\ref{eq:met}) on $\M$ is conformally flat.
\item[3).]
Equation (\ref{eq:ricl}) is satisfied everywhere on $\Pen$ if and only
if the traceless part $s_{ij}=r_{ij}-\frac{1}{4}g_{ij}r$ of the Ricci tensor 
of the metric (\ref{eq:met}) vanishes on $\M$. 
\end{itemize}
\et
A straightforward corollary from this theorem reads.
\bc~\\
Equation (\ref{eq:ricl}) is satisfied in $\Pen$ if and only if 
the base metric satisfies the Einstein equations $r_{ij}=\kappa g_{ij}$.
\ec
To interpret equation (\ref{eq:weyll}) geometrically on $\Pen$ it is 
convenient to consider it together with equations (\ref{eq:wc}) and 
(\ref{eq:wc'}). It is easily seen then that system 
(\ref{eq:wc})-(\ref{eq:weyll}) constitutes the Froebenius condition for 
the 3-dimensional distribution 
$\cal N$ which in $\Pen$ annihilates forms $\Lambda_L$, $F_L$ and $E_L$. It 
follows that $\cal N$ is totally null in the metric 
$\tilde{g}$ and has maximal dimension.\\ 

\noindent
We failed in finding a nice geometrical interpretation for 
equation (\ref{eq:ricl}). Since it is invariant under 
transformations (\ref{eq:trla})-(\ref{eq:trf}) such an interpretation should 
exist.

\section{Hermitian and optical geometries}

\noindent
Suppose that we are given 
a $2m$-dimensional real manifold $\cal R$ equipped 
with a real-valued metric $g$ 
of signature ($2p+\epsilon$, $2q+\epsilon$). 
Here $2m=2(p+q+\epsilon)$ and $\epsilon = 0$ or 1. Following 
\cite{bi:KT} we call the cases $\epsilon=0$ and $\epsilon=1$ 
pseudo-Euclidean and pseudo-Lorentzian, respectively. 
We omit the prefix ``pseudo'' if $pq=0$. 
Complexifying $g$ one endows 
the complexification  T${\cal R}^{\bf C}$ 
of the tangent bundle 
T$\cal R$ with a metric $g^{\bf C}$.
Let $\cal N$ be a vector subbundle of T$\cal R^{\bf C}$ which is totally 
null with respect to $g^{\bf C}$ and has $m$-dimensional fibers.  
We call such bundles maximal totally null bundles. 
Given $\cal N$ we also have its complex conjugate bundle $\bar{\cal N}$ 
as well as bundles 
${\cal N}\cap\bar{{\cal N}}$ and ${\cal N}+\bar{{\cal N}}$. 
It is easy to see that ${\cal N}\cap\bar{{\cal N}}$ and 
${\cal N}+\bar{{\cal N}}$ are respectively complexifications of certain  
vector subbundles $\cal K$ and ${\cal L}={\cal K}^\perp$ 
of T$\cal R$. The complex fiber dimension $r$ of 
${\cal N}\cap\bar{{\cal N}}$ 
(or real fiber dimension of $\cal K$) 
depends on the signature of $g$ and may take the following values 
$r=\epsilon,~2+\epsilon, ~..., 
\mbox{min}(2p+\epsilon,~2q+\epsilon)$ \cite{bi:KT}.  
It is 
called a real index of $\cal N$. From now on we only  
consider such $\cal N$s for which $r$ is constant over $\cal R$.\\
Given $\cal N$, we have a natural almost complex structure 
$\cal J$ in a bundle 
${\cal H}\stackrel{def}{=}{\cal L}/{\cal K}$. To define this we observe that
any section $l$ of $\cal L$ is of the form $l=n+\bar{n}$ where $n$ is some
section of $\cal N$. If $[l]$ denotes the equivalence class associated with
$l$ in $\cal H$ we define $\cal J$ by
$$
{\cal J}(~[l]~)={\cal J}(~[n+\bar{n}]~)\stackrel{def}{=}[-i~(n-\bar{n})~].
$$
One may prove that $\cal J$ is well defined. Moreover, since the metric $g$
is degenerate on $\cal K$ then it descends
to a unique metric $g'$ in $\cal H$. It further follows that 
$\cal J$ is an orthogonal 
transformation for $g'$ (we say that $\cal J$ is 
orthogonal with respect to $g'$ or, simply, orthogonal).\\
If $q=0$ then the only possible values of the real index $r$ of $\cal N$ are 0 
(Euclidean case) or 1 (Lorentzian case).
For such $q$ and $\epsilon =0$ we see that
the corresponding ${\cal K}=\{0\}$, ${\cal L}={\cal H}=$T${\cal R}$, 
$g'=g$. Thus, in this case, $\cal N$ defines an almost hermitian geometry
$(g, ~{\cal J})$ in $\cal R$.\\
If $q=0$ and 
$\epsilon =1$ then the maximal totally null bundle $\cal N$ equips 
$\cal R$ with the 
structure of the almost optical geometry of A. Trautman \cite{bi:ATOG,bi:T2}. 
This is a sequence
$$
\begin{array}{cccccc}
                          &{\cal K}&\hookrightarrow &{\cal L}&\longrightarrow 
&{\cal H}\\
\mbox{\rm fiber dimension}&1       &                &2m-1    &               
&2m-2
\end{array}
$$
of real vector subbundles $\cal K$ and $\cal L$ of T$\cal R$ together 
with an orthogonal 
almost complex structure $\cal J$ in the quotient bundle $\cal H$.
Note that in this case the metric $g'$ in $\cal H$ is purely Euclidean.\\
To deal with the generic case of $q\neq 0$ we introduce the following 
definition \cite{bi:NT}.
\bd~~\\
Consider a real $2m$-dimensional manifold $\cal R$ equipped with
a metric $g$ of signature $(2p+\epsilon,2q+\epsilon)$. 
Let $\cal K$ and ${\cal L}={\cal K}^{\perp}$ (${\cal K}\subset {\cal L}$) be 
vector subbundles of
${\rm T}\cal R$ which have respective fiber dimension $r$ and $2m-r$.
If the quotient bundle ${\cal H}={\cal L}/\cal K$ is equipped with an
almost complex structure $\cal J$ which is orthogonal with respect
to the descended metric $g'$ in $\cal H$, then
$({\cal K},g,{\cal H},{\cal J})$ is called an almost 
optical geometry with index $r$ (or 
almost $r$-optical geometry).
\ed

\noindent
Thus any maximal totally null bundle with real index $r$ gives rise
to an almost $r$-optical geometry. The converse is also true.
Given an almost $r$-optical geometry on $\cal R$ we define $\cal N$
as the vector subbundle of T${\cal R}^{\bf C}$ whose sections
are of the form $n=l_1+il_2,$ where $l_1,l_2$ are sections of $\cal L$
satisfying ${\cal J}[l_1]+i{\cal J}[l_2]=-i[l_1]+[l_2]$. One easily proves that
$\cal N$ so defined is totally null, has fiber dimension $m$ and real
index $r$.\\

\noindent
The above discussion shows a one-to-one correspondence between maximal totally
null bundles of a given index $r$ and almost $r$-optical geometries.\\

\noindent
Given a maximal totally null bundle $\cal N$ we denote the set of all its
sections by 
$\G ({\cal N})$. It is natural to consider the following integrability 
conditions for $\cal N$
\be
[\G ({\cal N}), \G ({\cal N})]\subset\G ({\cal N}).\label{eq:integr}
\ee
Here $[\cdot ,\cdot ]$ denotes a 
commutator of sections treated as vector fields. We say that an 
almost $r$-optical geometry 
associated with $\cal N$ is $r$-optical if and only if
the conditions (\ref{eq:integr}) are satisfied.\\ 

\bd~\\
A CR-structure \cite{bi:CR} is a real $(2m-1)-$dimensional manifold 
$\cal Q$ equipped with a subbundle $\cal H$ 
of the tangent bundle T$\cal Q$, which has fibers of 
dimension $2(m-1)$ and which is endowed with an almost complex 
structure $\cal J$. 
\ed
Given a CR-structure $({\cal Q}, {\cal H}, {\cal J})$ we extend 
$\cal J$ to the complexification ${\cal H}^\C$ by linearity. 
A CR-structure is called an integrable CR-structure if for any sections 
$X$, $Y$ of $\cal H$ we have 
\be
{\cal J}~~[X+i{\cal J}X,~Y+i{\cal J}Y] =
-i ~[X+i{\cal J}X,~Y+i{\cal J'}Y].\label{eq:cccr}
\ee

We say that two CR-structures $(\crs , ~{\cal H}, ~{\cal J})$ and
$(\crs ',~{\cal H'}, ~{\cal J'})$ are
(locally) equivalent
iff there exists a (local) diffeomorphism
$\phi:\crs\rightarrow\crs '$ such
that $$\phi _{*}{\cal H}={\cal H'}$$ and
$$\phi^*{\cal J}={\cal J'}.$$
In the following we will also need the more general 
structure. 
\bd~\\
An r-CR-structure is a real $(2m-r)$-dimensional manifold equipped
with a subbundle $\cal H$ of the tangent bundle T$\cal Q$ such that 
it has fibers of 
dimension $2m-2r$ and is endowed with an almost complex structure $\cal J$.\\
A $r$-CR-structure is integrable iff any two sections $X$, $Y$ of the bundle 
$\cal H$  
satisfy conditions (\ref{eq:cccr})
\ed
We note that 0-CR-structure is the same as an almost complex geometry in 
$\cal Q$. Its integrability conditions are equivalent to the integrability 
conditions of this 
almost complex structure.\\
Given an almost 
$r$-optical geometry $\og$ on $\cal R$ we choose a surface $\cal S$ 
of dimension $2m-r$ in $\cal R$ that it is 
transversal to sections of the bundle $\cal K$. It is easy to see that any 
such surface is naturally endowed with a 
$r$-CR-structure. 
If it happens that $\og$ is $r$-optical then 
the integrability conditions 
(\ref{eq:integr}) imply that $r$-CR-structures on  
any hypersurface $\cal S$ are integrable and locally equivalent. 
More formally, given an $r$-optical geometry satisfying (\ref{eq:integr}) 
we find that the bundle $\cal K$ is integrable as a distribution on $\cal R$. 
Thus it defines 
a foliation of $\cal R$ by $r$-dimensional real manifolds tangent to 
$\cal K$. Consider an equivalence relation $\sim$ in $\cal R$  
which identifies points lying on the 
same leaf $\cal X$ of this foliation. We assume that its quotient 
space $\crs\equiv{\cal R}/\sim$ is a manifold. Conditions (\ref{eq:integr}) 
guarantee that the projection of $r$-CR-structures from any of the surfaces 
$\cal S$ to this manifold equip it with the same  
integrable $r$-CR-structure. Hence, in such a case, the manifold 
$\cal R$ is locally equivalent to the Cartesian product  
${\cal X}\times {\cal Q}$. This generalizes the well known fact for almost 
optical geometries associated with congruences of shear-free 
and null geodesics in 4-dimensions \cite{bi:RT1}, \cite{bi:RT}.\\

\noindent
Summing up we have the following theorem.
\bt~\\
Let $\cal R$ be a real $2m$-dimensional manifold equipped with a real 
metric $g$ of signature ($2p+\epsilon$, $2q+\epsilon$), 
where $2m=2(p+q+\epsilon)$ and $\epsilon = 0$ or 1.
\begin{itemize}
\item[1)] There exists a one-to-one correspondence between
almost 
$r$-optical geometries over $\cal R$ and maximal totally null bundles 
${\cal N}$ of constant real index $r$ over ${\cal R}$.
\item[2)] Any integrable $\cal N$ of index $r$ 
locally defines an integrable $r$-CR-structure.
\item[3)] In the case of a Euclidean metric, the bundle 
$\cal N$ corresponds to an almost hermitian 
structure ($g,~{\cal J})$ on $\cal R$. This almost hermitian structure is 
integrable iff $\cal N$ satisfies integrability conditions 
(\ref{eq:integr}).
\item[4)] In the case of a Lorentzian metric, the bundle 
$\cal N$ corresponds to an almost optical 
geometry on $\cal R$. This, when integrable, defines an integrable 
CR-structure.
\end{itemize}
\et

\noindent
In the following we will need the interpretation of the integrability 
conditions
(\ref{eq:integr}) in terms of the theory of differential ideals.\\
Consider a system of complex-valued 1-forms $(A_1, A_2, ..., A_s)$ on
$\cal R$. Let $\cal I$ be an ideal in the exterior algebra of all
complex-valued differential forms on $\cal R$ generated by 1-forms
$(A_1, A_2, ..., A_s)$. We say that $\cal I$ is a closed differential ideal
iff
\beq
\der A_1\dz A_1\dz A_2\dz...\dz A_s=0,\nonumber\\
\der A_2\dz A_1\dz A_2\dz...\dz A_s=0,\nonumber\\
...\nonumber\\
...\nonumber\\
...\nonumber\\
\der A_s\dz A_1\dz A_2\dz...\dz A_s=0.\nonumber
\eeq

\noindent
Any maximal totally null bundle $\cal N$ over $\cal R$ 
can be defined as the  annihilator of 
$m$ linearly independent, totally null, complex-valued 1-forms, say  
$(A_1, A_2, ..., A_m)$, on $\cal R$. Given $\cal N$ defined by such 1-forms 
we have the following, well known, lemma.
\bl~\\
$\cal N$ satisfies the integrability conditions (\ref{eq:integr}) if and 
only if  
the system $(A_1, A_2, ..., A_m)$ generates a closed differential ideal 
on $\cal R$.
\el
\noindent

\section{Weyl geometry.}

\noindent
{\bf 4.1 Definitions}\\

\noindent
From now on by a metric on a real manifold we will understand 
a nondegenerate, bilinear and symmetric, complex-valued form. \\

\noindent
Consider a 4-dimensional real oriented manifold $\M$ equipped 
with a metric $g$. Fixing four complex-valued 1-forms 
$(\theta^i)$ $(i=1,2,3,4)$ on $\M$ for which 
$\theta^1\dz\theta^2\dz\theta^3\dz\theta^4\neq 0$ we can represent any metric 
$g$ 
by means of its 
metric coefficients $g_{ij}$. Thus, given $g$ and $(\theta^i)$ we have 
$$g=g_{ij}\theta^i\theta^j.$$ 
The system of forms $(\theta^i)$ will be called a cotetrad on $\M$. 
We equip $\M$ with a Weyl geometry. Such geometry is defined in terms of 
a pair $(g, A)$ where $g$ is a metric and $A=A_i\theta^i$ is a complex-valued 
1-form on $\M$. The metric and $A$ are related by 
\be
Dg_{ij}=dg_{ij}-g_{ik}\Gamma^k_{~j}-g_{jk}\Gamma^k_{~i}
=-Ag_{ij}\label{eq:weylg},
\ee 
where $\Gamma^i_{~j}$ are torsion-free connection 1-forms. The torsion-free 
condition is expressed by   
\be
\der\theta^i=-\Gamma^i_{~j}\dz\theta^j.
\label{eq:skr}
\ee
Given a Weyl geometry $(g,A)$ on $\M$ the connection 1-forms $\G^i_{~j}$ are 
uniquely determined. They are expressible in terms of $A$ and the Levi-Civita 
connection 1-forms $\g^i_{~j}$ of the metric $g=g_{ij}\theta^i\theta^j$. 
Explicitely we have 
\be
\Gamma^i_{~j}=g^{ik}\Gamma_{kj}, ~~~~~~~~g^{ik}g_{kj}=\delta^i_{~j}
\label{eq:wlc1}
\ee
where
\be
\Gamma_{ij}=\gamma_{ij}+\frac{1}{2}g_{ij}A+ 
g_{k[i}A_{j]}\theta^k,~~~~~~~\g_{ij}=g_{ik}\g^k_{~j}
\label{eq:wlc}
\ee
and where we have introduced the abbreviation 
$a_{[i}b_{j]}=\frac{1}{2}(a_ib_j-a_jb_i)$
\footnote{Round bracket will denote symmetrization of 
indices, e.g. $a_{(i}b_{j)}=\frac{1}{2}(a_ib_j+a_jb_i)$.}.

\noindent
Given connection 1-forms $\Gamma^i_{~j}$ one may associate a (Weyl) 
connection with them and obtain a receipee for parallel transport of 
vectors on $\M$. It follows that in
contrast to the parallel transport of Riemannian geometry, this transport 
preserves only nullity of vectors 
(see, for instance, \cite{bi:Weyl} for more information)
\footnote{Our point of view on Weyl geometries is nonstandard in two respects. 
Firstly, we 
admit complex metrics $g$. Secondly, we do not stress the conformal 
invariance. It is easy to check that for a fixed cotetrad $(\theta^i)$ 
equation (\ref{eq:weylg}) is invariant under the transformation 
$$(g, A)\rightarrow (g', A')=
({\rm e}^{2\varphi}g,~ A-2{\rm d}\varphi).~~~~~~~~~~~~~~~~~~~~~~~(*)$$
One can therefore view Weyl 
geometry as a pair ($g, A$) given up to transformations $(*)$. In such a 
formulation only a conformal metric is relevant. We do not refer to this 
point of view 
in our discussion since we want to have a nice passage to Riemannian 
geometries (fixed metrics, not their conformal class) when $A=0$. 
However, all relevant formulae such as for example 
(\ref{eq:eiwey}), (\ref{eq:integr1}), 
(\ref{eq:integr2}), (\ref{eq:weyl}) and (\ref{eq:ric}) are covariant under 
$(*)$. See Appendix E for a further discussion of this.
\\Thus our results apply also to Weyl geometries viewed in this
standard, conformal sense.}.

\noindent
The curvature of Weyl geometry is defined in terms of 
curvature 2-forms 
\be
\Omega_{ij}=
\frac{1}{2}R_{ijkl}\theta^k\dz\theta^l=
\der\G_{ij}+\G_{ik}\dz\G^k_{~j}.
\ee
It splits into 
the curvature $\omega_{ij}$ of the Levi-Civita connection, and 
the remaining $A$-dependent part. This, in particular, includes the 
curvature 
$${\cal F}=\frac{1}{2}{\cal F}_{ij}\theta^i\dz\theta^j=\der A$$ 
of $A$.\footnote{Decompositions of various Weyl geometric objects 
onto the Levi-Civita and $A$-dependent parts are given in Appendix A.}\\
The Ricci tensor $R_{jl}$ and its scalar $R$ is defined respectively by 
$R_{jl}=g^{ik}R_{ijkl}$ and $R=g^{ij}R_{ij}$. Note that 
$R_{ij}$ is not symmetric in general. The 
traceless part of $R_{(jl)}$ is defined by 
\be
S_{ij}=R_{(ij)}-\frac{1}{4}Rg_{ij},
\ee 
which implies
\be
S:=g^{ij}S_{ij}=0.\label{eq:slad}
\ee
We say that Weyl geometry $(g,A)$ satisfies the Einstein-Weyl equations iff 
\be
S_{ij}=0.\label{eq:eiwey}
\ee
For further use we also define a tensor 
\be
C_{ijkl}=R_{ijkl}+\frac{1}{3}Rg_{i[k}g_{l]j}+R_{j[k}g_{l]i}+R_{i[l}g_{k]j}.
\ee 
This can be decomposed into the Levi-Civita 
($w_{ijkl}$) and $A$-dependent ($f_{ijkl}$) 
parts  
\be
C_{ijkl}=w_{ijkl}+f_{ijkl}
\ee
(see Appendix A). It turns out that $w_{ijkl}$ are precisely the covariant 
coefficients of the usual 
Weyl tensor associated with the metric $g$. 
Since $w_{ijkl}$ are antisymmetric in $k,l$ then we may associate with them 
a collection of 2-forms 
\be
w_{ij}=\frac{1}{2}w_{ijkl}\theta^k\dz\theta^l,
\ee
which we call the Weyl-tensor 2-forms.\\

\noindent 
We close this section with a remark that if $A=0$ everywhere on $\M$ then 
such Weyl geometry reduces to the usual Riemannian geometry associated with 
metric $g$. In particular, such objects like $\G_{ij}$, $R_{ijkl}$ etc. 
reduce to their respective Levi-Civita parts $\g_{ij}$, $r_{ijkl}$ etc. 
(compare with Footnote 2). \\

\noindent
{\bf 4.2 Weyl geometries in null tetrads}\\

\noindent
Of particular 
interest are null cotetrads on $\M$. These are cotetrads
\be
(\theta^1, ~\theta^2, ~\theta^3, \theta^4)=(M,~P,~N,~K)
\label{eq:nt}
\ee 
related to the metric $g$ by  
\be
g=g_{ij}\theta^i\theta^j=MP-NK.\label{eq:met1}
\ee
A tetrad dual to (\ref{eq:nt}) will be denoted by 
\be
(e_1,~e_2,~e_3,~e_4)=(m,~p,~n,~k).\label{eq:ntv}
\ee

From now on we restrict only to null cotetrads that agree with the orientation 
i.e. such null cotetrads for which volume form $\eta$ of $\M$ reads
\be
\eta=M\dz P\dz N\dz K.\label{eq:volume}
\ee 
 
Given a null cotetrad and $A$ we find Weyl connection 1-forms $\G^i_{~j}$, and 
calculate curvature 2-forms. Their convenient decomposition 
relates to the notion of selfduality. \\

\noindent
Given a $p$-form $\omega$ on $\M$ we define its Hodge 
dualization $*\omega$ by  
\be
(*\omega)(X_1,...,X_{4-p})\eta=
\omega\dz g(X_1)\dz ...\dz g(X_{4-p}), 
\ee
where $g(X_i)$ is a 1-form associated with a vector field $X_i$ 
$(i=1,2,...,(4-p))$ by $<g(X_i),X_j>=g(X_i,X_j)$. Since 
the metric $g$ induces an isomorphism between forms and vectors on $\M$ 
then, in an obvious way, we have also a Hodge dualization of $p$-vectors.

\noindent
Hodge dualization is an involutive ($*^2=$id) automorphism of the complexified 
space $\bigwedge^2$ of 2-forms on $\M$. Its $\pm$ eigenspaces 
$\bigwedge^2_+$ and $\bigwedge^2_-$ consist of self-dual and 
antiself-dual forms respectively. A convenient basis for $\bigwedge^2_+$ is 
\be
P\dz K,~~~~~N\dz K-M\dz P,~~~~~N\dz M
\ee
and for $\bigwedge^2_-$ is
\be
M\dz K,~~~~~N\dz K+M\dz P,~~~~~N\dz P.
\ee
Any 2-form can be decomposed onto these bases. Decompositions of 
the curvature $\cal F$ and 
the Weyl tensor 2-forms $w_{ij}$ onto these bases define 
coefficients 
$\phi_0$, $\phi_1$, $\phi_2$, $\phi'_0$, $\phi'_1$, $\phi'_2$  
and $\Psi_0$, $\Psi_1$, $\Psi_2$, $\Psi_3$, $\Psi_4$, 
$\Psi'_0$, $\Psi'_1$, $\Psi'_2$, $\Psi'_3$, $\Psi'_4$ by  
\be
{\cal F}=\phi'_0 N\dz P+\phi'_1(N\dz K+M\dz P)+\phi'_2 M\dz K+\phi_0 N\dz M 
+\phi_1 (N\dz K - M\dz P)+\phi_2 P\dz K,
\ee
\beq
w_{14}&=\Psi'_0 M\dz K+\Psi'_1 (N\dz K+M\dz P)
+\Psi'_2N\dz P\nonumber\\
w_{23}&=-\Psi'_2 M\dz K-\Psi'_3 (N\dz K+M\dz P)
-\Psi'_4N\dz P\label{eq:wey2+}\\
\frac{1}{2}(w_{34}+w_{12})&=\Psi'_1 M\dz K+\Psi'_2 (N\dz K+M\dz P)+
\Psi'_3 N\dz P\nonumber
\eeq
\beq
w_{24}&=\Psi_0 P\dz K+\Psi_1 (N\dz K-M\dz P)+
\Psi_2 N\dz M\nonumber\\
w_{13}&=-\Psi_2 P\dz K-\Psi_3 (N\dz K-M\dz P)-
\Psi_4 N\dz M\label{eq:wey2-}\\
\frac{1}{2}(w_{34}-w_{12})&=\Psi_1 P\dz K+\Psi_2 
(N\dz K-M\dz P)+\Psi_3 N\dz M\nonumber. 
\eeq
Decompositions of $\Omega_{[ij]}$ can be  
found in Appendix B.\\
It follows from the above decompositions 
that the Weyl tensor 2-forms $w_{ij}$ are anti-self-dual (respectively, 
self-dual) if and only if all the coefficients 
$\Psi_0$, $\Psi_1$, $\Psi_2$, $\Psi_3$, $\Psi_4$ 
(respectively, $\Psi'_0$, $\Psi'_1$, $\Psi'_2$, $\Psi'_3$, $\Psi'_4$) vanish
\footnote{It is known that conditions $\Psi_0=\Psi_1=\Psi_2=\Psi_3=\Psi_4=0$ 
or $\Psi'_0=\Psi'_1=\Psi'_2=\Psi'_3=\Psi'_4=0$ are 
invariant under the conformal transformations of the metric. A less well 
known fact 
states that they are also invariant under transformations $(*)$ of 
Footnote 2. This is related to the fact that the conditions $C_{ijkl}=0$ and 
$w_{ijkl}=0$ (hence also  
$f_{ijkl}=0$) are invariant under $(*)$.} on 
$\M$. This allows for the following terminology. Weyl geometries $(g, A)$ 
on $\M$ are 
called anti-self-dual (respectively, self-dual) if and only if all 
coefficients $\Psi_\mu$ (respectively $\Psi'_\mu$), $\mu=0,1,2,3,4$ vanish on 
$\M$. It follows that this definition does not depend on the choice of 
the null cotetrad.

\section{Geometry of twistor bundles}
{\bf 5.1 Twistor bundles}\\

\noindent
Let $\M$ be a real oriented 4-manifold equipped with a Weyl geometry 
$(g,A)$. At any point $x$ of $\M$ we consider vector subspaces of the 
complexification of the tangent space T$_x\M$ which
\begin{itemize}
\item[1)] are totally null with respect to $g$,
\item[2)] and have maximal dimension.
\end{itemize}
Such spaces are necessarily 2-dimensional and can be represented by a 
complex bivector. It turns out that bivectors associated with spaces 
satisfying 1) and 2) are either self-dual or anti-selfdual in the 
Hodge dualization associated with $g$ and the orientation of $\M$. This 
shows that 
the set of all spaces that at $x$ satisfy 1) and 2) consists of two disjoint 
parts 
$S_x$ and $S'_x$. We call $S_x$ (respectively, $S'_x$) a set of all 
self-dual (anti-self-dual) maximal totally null spaces at $x$. A pair 
of spaces $(s,s')$ such that $s\in S_x$ and $s'\in S'_x$ is called 
a pair of maximal totally null spaces of opposite self-duality. 
It is easy to see that both $S_x$ and 
$S'_x$ are diffeomorphic to a 2-dimensional sphere ${\bf S}_2$. 
A stereographic projection gives a convenient parametrization of these 
spheres in terms of points of the Argand plane ${\bf C}\cup\{\infty\}$. 
Using the null 
tetrad (\ref{eq:ntv}) for $g$ we find that elements $s\in S_x$ and 
$s'\in S'_x$ can be respectively 
represented by
\be
s={\rm Span}\{zm-k,zn-p\}~~~~~~~s'={\rm Span}\{z'p-k,z'n-m\},
\label{eq:ss'}
\ee
where $z,z'\in{\bf C}\cup\{\infty\}$ and thus we identified points of 
${\bf S}_2$ 
with the points of ${\bf C}\cup\{\infty\}$. For further use we also note that 
any two spaces $s$ and $s'$ have nonzero intersection at $x$. It is easy to 
see that this intersection is 1-dimensional and is spanned by a null 
vector $X$ which, if $s$ and $s'$ are represented by (\ref{eq:ss'}), has the  
form 
\be
X=k+zz'n-zm-z'p.\label{eq:X}
\ee

\noindent
Collecting 
the sets $S_x$, $S'_x$ point by point we have two fiber bundles 
${\cal P}=\displaystyle{\bigcup_{x\in\cal M}}S_x$ and 
${\cal P}'=\displaystyle{\bigcup_{x\in\cal M}}S'_x$ over $\M$. Any of these 
bundles has ${\bf S}_2$ as its typical fiber and is 
equipped with respective projections 
$\pi:\Pen\rightarrow\M$ and $\pi':\Pen'\rightarrow\M$. Any point $p$ of $\Pen$ 
is a certain totally null (necessarily self-dual) space of maximal dimension 
at the corresponding point $x$ of $\M$ (there is an analogous statement for 
points of $\Pen'$). It can therefore be parametrized by $(x,z,\bar{z})$, where 
$z$ is as in (\ref{eq:ss'}). A point $p'\in\Pen'$ is parametrized 
by $(x,z',\bar{z}')$, where $\pi'(p')=x$ and $z'$ is as in (\ref{eq:ss'}). 
Using 
$\Pen$ and $\Pen'$ 
one defines their fibre product 
$\pp=\displaystyle{\bigcup_{x\in\cal M}}(S_x\times S'_x)$ which is 
\begin{itemize}
\item[1)] a fiber bundle over $\M$ with a natural projection  
$\Pi:\pp\rightarrow \M$ and typical fiber diffeomorphic to 
${\bf S}_2\times{\bf S}_2$,
\item[2)] a fiber bundle over $\Pen$ with a natural projection  
${\rm pr}:\pp\rightarrow \Pen$ and typical fiber diffeomorphic to 
${\bf S}_2$,
\item[3)] a fiber bundle over $\Pen'$ with a natural projection  
${\rm pr}':\pp\rightarrow \Pen'$ and typical fiber diffeomorphic to 
${\bf S}_2$.
\end{itemize}
In particular, any point of $\pp$ may be understood as 
a pair of maximal totally null 
spaces of opposite self-duality at the corresponding point of $\M$.
A convenient parametrization of $\pp$ is $(x,z,\bar{z},z',\bar{z}')$.
\noindent
The projections associated with these bundles satisfy 
\be
\Pi=\pi\circ{\rm pr}=\pi'\circ{\rm pr}',\label{eq:pull}
\ee 
and in the above coordinates on $\pp$, $\Pen$ and $\Pen'$ are given by 
$$\Pi(x,z,\bar{z},z',\bar{z}')=x,~~~~~~~
{\rm pr}(x,z,\bar{z},z',\bar{z}')=(x,z,\bar{z}),~~~{\rm etc.}$$
Using the projections we can pull back forms. For example, 
using $\pi$ we pull back forms from $\M$ to $\Pen$, using $\Pi$ we 
pull back forms from $\M$ to $\pp$ and using pr$'$ we pull back 
forms from $\Pen'$ to $\pp$. \footnote{Note that due to property 
(\ref{eq:pull}) the direct 
pull back of a form from $\M$ to $\pp$ is the same as a pull back of 
this form via $\Pen$ (first using $\pi$ and then pr) or $\Pen'$.} 
In this way we can for example pullback tetrad 1-forms 
$\theta^i$ and Weyl connection 1-forms $\Gamma^i_{~j}$ from $\M$ to 
$\Pen$, $\Pen'$ and $\pp$. Since it follows from the 
context on which manifold a given form is placed, we use in the 
following the same letters to denote forms and 
their pullbacks.

\noindent
Weyl geometry $(g, A)$ induces interesting geometrical structures 
on bundles $\Pen$, $\Pen'$ and $\pp$. We only outline constructions for $\Pen$ 
and $\pp$. \\

\noindent
{\bf 5.2 Natural structures on $\Pen$ and $\Pen'$}\\

\noindent
{\it 1) The tangent bundle to $\Pen$ and its complexification 
split naturally into vertical and horizontal part.} \\

\noindent
To show this we give a receipee for 
the lifting of a given vector $v$ from $x\in\M$ to a chosen point $p\in\Pen$ 
in the fiber over $x$. Recall that a point $p$ can be considered a 
maximal totally null self-dual space at $x$. Take 
any curve $x(t)$ that is tangent to $v$ at $x$. Use the Weyl connection 
associated with $(g, A)$ to propagate the maximal totally null space 
represented by $p$ paralelly 
along $x(t)$. Since Weyl-geometric 
parallel propagation preserves nullity of vectors, 
then at  
any point of our curve we get a certain 
totally null space. Due to the continuity of $x(t)$ any such space is 
maximal and self-dual. Thus given a 
curve $x(t)$ tangent to $v$ at $x$ we have a corresponding curve $p(t)$ 
in $\Pen$ which starts at $p$. It follows that a 
direction of a 
tangent vector to $p(t)$ in $p$ does not depend on the choice of $x(t)$. 
This is the direction of the desired 
horizontal lift $\tilde{v}$ of $v$ to $p$. The lift is determined 
completely by the additional demand that $\pi_*(\tilde{v})=v$. 
Thus we are able to lift any vector from $x\in\M$ to 
a chosen point $p$ from the fiber over $x$. It is moreover true that 
horizontal lifts of four linearly independent vectors from $x$ constitute 
four linearly independent vectors in $p$. This means that we have a well 
defined 
lift of the tangent space T$_x\M$ to a 4-dimensional subspace 
$H_p$ of T$_p\Pen$. This subspace is called a 
horizontal space at $p$. The vertical space $V_p$ consists of 
vectors at $p$ that are tangent to the fibers. This space is 2-dimensional and 
may be identified with the tangent space to a certain point of ${\bf S}^2$. 
A direct sum $H_p\oplus V_p$ equals to T$_p\Pen$.

\noindent
The horizontal lift that we described above can be also used to lift 
horizontally vectors $w$ from the complexification of the tangent space 
at $x$ to $p\in\pi^{-1}(x)\subset\Pen$. This is achieved by extending  
the horizontal lift map $\tilde{~}: v\rightarrow\tilde{v}$ by linearity to the 
complexification of T$_x\M$. Thus $\tilde{w}=\tilde{v}_1+i\tilde{v}_2$, where 
$v_1$ and $v_2$ are respectively real and imaginary parts of $w$. This 
leads to a split of the complexification of the tangent bundle to $\Pen$ 
into vertical and horizontal parts 
$(\rm{T}\Pen)^{\bf C}={\cal H}^{\bf C}\oplus{\cal V}^{\bf C}$.

\noindent
There is a natural complex structure $I_p$ 
in $V_p$ that comes from the natural complex  
structure on ${\bf S}^2$. This, when prolonged to the complexification 
$V^{\bf C}_p$, gives a split $V^{\bf C}_p=V^+_p\oplus V^-_p$, where    
$I_pV^{\pm}_p=\pm iV^{\pm}_p$.\\

\noindent
{\it 2) Spin connection 1-form.}\\

\noindent
We look for a complex-valued 1-form $E$ 
such that in some neighbourhood ${\cal U}$ of $\cal P$ it satisfies 
$i)$ $E({\cal H})=0$, $ii)$ $E(\displaystyle{\bigcup_{p\in\cal U}}V^-_p)=0$ 
and 
$iii)$ $E\dz \bar{E}\neq 0$. In general, starting from a given point 
$x\in \cal M$, we can solve these conditions only in 
a cylindrical ${\cal U}$ over a sufficiently small 
neighbourhood of $x$. Outside $\cal U$ conditions 
$i)-iii)$ may be contradictory. From now on we restrict our 
considerations to such $\cal M$s for which the corresponding 
$E$ is defined globally. This can be achieved either by assuming 
some additional structure on $\cal M$ (e.g. spin manifold structure) 
or restricting to $\cal M$s to be open subsets of ${\bf R}^4$. With  
such an assumption, conditions $i)-iii)$ define $E$ on $\Pen$ 
up to a multiplication by a 
non vanishing complex-valued function $h$ on $\Pen$ 
\be
E\rightarrow hE.\label{eq:e}
\ee
$E$ is called a spin connection 1-form on $\Pen$.\\
Using the null cotetrad (\ref{eq:nt}) and the coordinates $(x,z,\bar{z})$ on 
$\Pen$
one easily finds that the form $E$ may be represented by 
\be
E=dz-\Gamma^3_{~2} 
+z(\Gamma^1_{~1}-\Gamma^4_{~4})+z^2\Gamma^2_{~3}. 
\ee

\noindent
{\it 3) Metrics.}\\

\noindent
Pullback the metric $g$ from $\M$ to $\Pen$ and add to it a tensor 
$h\bar{h}E\bar{E}$ with $h$ being a non vanishing function on $\Pen$. This 
defines a class of metrics $\tilde{g}$ on $\Pen$, 
which can be represented by  
\be
\tilde{g}=\pi^*(g)+h\bar{h}E\bar{E}.\label{eq:gtil}
\ee\\

\noindent
{\it 4) Canonical field of 2-planes.}\\

\noindent
Take a point $p$ of $\Pen$. It corresponds to a certain self-dual maximal 
totally null plane at $x=\pi(p)\in \cal M$. Lift this plane horizontally to 
$p$. This associate a horizontal 2-plane to any point $p$ of $\Pen$. Thus on 
$\Pen$ we have a distinguished field of 2-planes,   
which we call the canonical field of 2-planes. Note that any 2-plane 
in this field is totally null in any metric from the class (\ref{eq:gtil}).

\noindent
Dually, the canonical field 
of 2-planes defines a pair of 1-forms $(F, B)$ on $\Pen$ which, 
by definition annihilate the vertical space and the canonical field of 
2-planes. 
$F$ and $B$ are given 
up to transformations
\be
(F,B)\rightarrow (\al F+\bet B, \g F+\delta B),~~~~\al\delta-\bet\g\neq 0.
\label{eq:fb}
\ee
This shows that a {\it direction} of a 
2-form $\Sigma=F\dz B$ on $\Pen$ is well defined.\\
It is easy to see that in the null cotetrad (\ref{eq:nt}) and 
in the coordinates $(x,z,\bar{z})$ the forms $F$ and $B$ may be represented by 
\be
F=M+zK, ~~~~~~B=N+zP.\label{eq:fbb}
\ee\\
\noindent
{\it 5) Distinguished totally null planes of maximal dimension.}\\

\noindent
Given a point $p\in\Pen$ consider a canonical 2-plane $\sigma$ passing 
through this point. 
There are only two 3-dimensional planes at $p$ that are totally null in any 
metric $\tilde{g}$ and that contain $\sigma$ as a vector 
subspace. These may be defined  
as vector spaces $n_E$ and $n_{\bar{E}}$ 
annihilating $(F, B, E)$ and $(F, B, \bar{E})$, respectively. 
Point by point they define two bundles of maximal 
totally null planes ${\cal N}_E$ and ${\cal N}_{\bar{E}}$ over $\Pen$. 
According to Section 2, in  
the case of real geometries, they will define a pair of distinguished 
complex or optical structures on $\cal P$.\\

\noindent
{\it 6) Invariant equations.}\\

\noindent
Although forms $E$, $F$, $B$ and $\Sigma$ are only given up to certain 
transformations, one can use them to write down several geometric 
equations on $\Pen$. In particular, note that the equation
\be
\der E\dz F\dz B\dz E=0
\ee
as well as the system of equations 
\be
\der F\dz F\dz B\dz E=0~~~~~~\der B\dz F\dz B\dz E=0
\ee
are invariant under the transformations (\ref{eq:e}), (\ref{eq:fb}).\\

\noindent
$1'), 2'), 3'), 4'), 5'), 6').$ Analogous constructions as in 1), 
2), 3), 4), 5), 6) can be 
performed for $\Pen'$. In this way on $\Pen'$ we have a split of T$\Pen'$ 
into vertical and horizontal part. Also the spin 
connection 1-form $E'$, metrics,  
canonical field of 2-planes, classes 
of forms $(F', B')$, $\Sigma'$, distinguished  
maximal totally null planes and 
invariant equations are defined there.\\

\noindent
{\bf 5.3 Natural structures on $\pp$.}\\

\noindent
{\it 1) The tangent bundle to $\pp$ and its complexification 
have a natural split into vertical and horizontal parts.} \\

\noindent
The receipee for having this split is almost like in the case of $\Pen$ with 
the exception that now a point $p\in\pp$ corresponds to a pair 
$(s,s')$ of maximal totally null 
spaces of opposite self-duality at $x=\Pi(p)\in\M$. Thus if we 
want to lift a vector $v$ from $x\in\M$ to a point $p\in\Pi^{-1}(x)\subset\pp$ 
we take a curve tangent to $v$ at $x$ and propagate paralelly spaces 
$s$ and $s'$ 
along this curve. This produces a pair 
of maximal totally null spaces of opposite self-duality at any point along 
the curve. Correspondingly we 
get a curve in $\pp$ starting at $p$ which defines 
the direction of the lift $\tilde{\tilde{v}}$ of $v$. As before the lift is 
specified 
uniquely by the demand that $\Pi_*(\tilde{\tilde{v}})=v$. 
Lifting T$_{\Pi(p)}\M$ horizontally we get a horizontal space $H_p$ in $p$. 
The vertical space $V_p$ is defined as the vector space tangent at $p$ to 
the fiber of $\pp$ over $\Pi(p)$. 
Note that now $V_p$ 
is 4-dimensional and is isomorphic to the tangent 
space of ${\bf S}^2\times{\bf S}^2$ at the point corresponding to $p$. \\

\noindent
{\it 2) Connection 1-forms.}\\

\noindent
These are the complex-valued 1-forms on $\pp$ that 
annihilate the horizontal space in T$(\pp)$. It follows that the basis of 
such forms 
on $\pp$ is given by the four pullbacks 
pr$^*(E)$, pr$^*(\bar{E})$, pr$'^*(E')$ 
and pr$'^*(\bar{E}')$. Here we wrote pullback signs explicitly. They will 
be omitted in the following.\\
Since the form $E$ (respectively $E'$)
was defined on $\Pen$ (respectively on $\Pen'$) up to a scaling by a 
function, 
the four above mentioned forms are given up to a scaling by a 
non vanishing complex function on $\pp$.\\
Local representations of $E$ and $E'$ may be given in terms of the coordinates 
$(x,z,\bar{z},z',\bar{z}')$ introduced on $\pp$ in Section 4.1. Since they 
were chosen in such a way that by projections 
we were getting corresponding coordinates $(x,z,\bar{z})$ on $\Pen$ 
and $(x,z',\bar{z}')$ on $\Pen'$ then we easily find that 
\be
E=dz-\Gamma^3_{~2} 
+z(\Gamma^1_{~1}-\Gamma^4_{~4})+z^2\Gamma^2_{~3},
\ee
\be
E'=dz'-\Gamma^3_{~1} 
+z'(\Gamma^2_{~2}-\Gamma^4_{~4})+z'^2\Gamma^1_{~3}.
\ee
Here, as usual, connection 1-forms are expressed with respect to the cotetrad 
(\ref{eq:nt}).\\

\noindent
{\it 3) Metrics.}\\

\noindent
The following metrics are of particular interest on $\pp$.
\be
\tilde{\tilde{g}}=\Pi^*(g)+h\bar{h}E\bar{E}+h'\bar{h}'E\bar{E}',
\label{eq:gttil}
\ee
where $h$ and $h'$ are non vanishing complex-valued functions on $\pp$.\\

\noindent
{\it 4) The canonical field of 3-planes and associated bundles.}\\

\noindent
At every point $p$ of $\pp$ there is a natural 3-plane which is obtained as 
follows. Take a pair $(s,s')$ of maximal totally null spaces of opposite 
self-duality which at $x=\Pi(p)$ correspond to $p$. Lift spaces 
$s$ and $s'$ horizontally to a point $p\in\Pi^{-1}(x)$ corresponding 
to $(s,s')$. This gives a pair of vector spaces 
$\tilde{\tilde{s}}$ and $\tilde{\tilde{s}}'$ at $p$. But as we noticed in 
Section 4.1 $s$ and $s'$ have a 1-dimensional intersection. Hence the vector 
space $\tilde{\tilde{s}}+\tilde{\tilde{s}}'$ has complex 
dimension equal to three. Thus at 
every point $p\in\pp$ we have a 3-dimensional space 
$\tilde{\tilde{s}}+\tilde{\tilde{s}}'$, 
which we call a canonical field of 3-planes.

\noindent
Actually the above considerations show that we have a list of vector 
subbundles ${\cal S}=\displaystyle{\bigcup_{p\in\pp }}\tilde{\tilde{s}}$, 
${\cal S'}=\displaystyle{\bigcup_{p\in\pp }}\tilde{\tilde{s}}'$, 
${\cal L}=\displaystyle{\bigcup_{p\in\pp }}(\tilde{\tilde{s}}
+\tilde{\tilde{s}'})$, 
${\cal K}=\displaystyle{\bigcup_{p\in\pp }}(\tilde{\tilde{s}}\cap 
\tilde{\tilde{s}}')$ 
of the complexification of the tangent bundle to $\pp$ which give rise 
to the following sequence.

$$
\begin{array}{cccccc}
                          &{\cal K}&\hookrightarrow &{\cal L}
&\longrightarrow &{\cal L}/\cal K\\
\mbox{\rm fiber dimension}&1       &                &3   &               
&2
\end{array}
$$

\noindent
Note that by definition bundles $\cal S$, ${\cal S}'$ are subbundles 
of $\cal L$, and that $\cal S$, ${\cal S}'$, and $\cal K$ are null 
bundles with respect to any metric from the class $\tilde{\tilde{g}}$. This 
indicates parallels between the structures defined here and the optical 
geometries of A. Trautman \cite{bi:ATOG,bi:T2}.

\noindent
Given the above bundles on $\pp$ it is interesting to ask whether such 
geometric conditions like $[{\cal S},{\cal S}]\subset\cal S$, 
$[{\cal S'},{\cal S'}]\subset\cal S'$, $[{\cal S},{\cal S}']\subset\cal L$, 
$[{\cal K},{\cal L}]\subset\cal K$ 
etc. have some interpretation in terms of the Weyl geometry on $\M$.\\

\noindent
{\it 5) Canonical 1-form on $\pp$.}\\

\noindent
The bundle $\cal K$ has 1-complex-dimensional fibers. It can be used to define 
a direction of 1-form $\Lambda$ on $\pp$. Indeed, if $X$ is a section 
of $\cal K$, then we define $\Lambda$ by 
\be
\Lambda=\tilde{\tilde{g}}(X).\label{eq:l}
\ee
Taking another section of $\cal K$ we see that 
\be
\Lambda\rightarrow u\Lambda, \label{eq:ltr}
\ee
which shows that a direction of $\Lambda$ is well defined. We call $\Lambda$ 
a canonical 1-form. We notice that $\La$, together with $E$, $\bar{E}$, 
$E'$, $\bar{E}'$, can be used to define a convenient basis of one 
forms on $\pp$. Indeed, one easily finds that 
\be
\tilde{\tilde{g}}=h\bar{h}E\bar{E}+h'\bar{h}'E'\bar{E}'+\Lambda T + 
FF',\label{eq:ff'}
\ee
with some complex valued 1-forms $F$, $F'$ and $T$ on $\pp$. These forms 
are defined up to the following transformations
\be
F\rightarrow\al F +\bet \Lambda,
\ee
\be
F'\rightarrow\frac{1}{\al}F'+\g\Lambda,
\ee
\be
T\rightarrow\frac{1}{u}(T-\al\g F-\frac{\bet}{\al}F'-\bet\g\Lambda ),
\ee
where $\al\neq 0$, $\bet, \g$ are some functions on $\pp$. It follows that 
the forms $F$ and $F'$ are in the class of forms obtained by taking pullbacks 
of the forms $F$ and $F'$ of Section 4.2. The set of forms 
$(E,\bar{E},E',\bar{E}',F,F',\Lambda, T)$ constitutes a convenient basis 
of 1-forms on $\pp$. In the coordinates $(x,z,\bar{z},z',\bar{z}')$ and in  
the cotetrad (\ref{eq:nt}) they may be represented as  
\beq
&E=dz-\Gamma^3_{~2} 
+z(\Gamma^1_{~1}-\Gamma^4_{~4})+z^2\Gamma^2_{~3},\nonumber\\
&E'=dz'-\Gamma^3_{~1} 
+z'(\Gamma^2_{~2}-\Gamma^4_{~4})+z'^2\Gamma^1_{~3},\nonumber\\
&\Lambda=-N-zz'K-zP-z'M,\nonumber\\
&F=M+zK,\label{eq:fonpp}\\
&F'=P+z'K,\nonumber\\
&T=K.\nonumber
\eeq\\
It is useful to consider transformations 
\beq
&z\leftrightarrow z',\label{eq:trpr}\\
&1\leftrightarrow 2,~~~~3\leftrightarrow 3, ~~~~4\leftrightarrow 4\nonumber
\eeq
where the last transformation means: for any object on $\cal M$, $\Pen$, 
$\Pen'$ or $\pp$, with tetrad indices 
1, 2, 3, 4 interchange indices 1 and 2 and do not change indices 3 and 4.\\
Examples: 
\beq
&\theta^1=M\leftrightarrow \theta^2=P, ~~~~~\theta^3=N\leftrightarrow 
\theta^3=N, ~~~~~\theta^4=K\leftrightarrow \theta^4=K,\\
&e_1=m\leftrightarrow e_2=p, ~~~~~e_3=n\leftrightarrow 
e_3=n, ~~~~~e_4=k\leftrightarrow e_4=k,\\
&S_{13}\leftrightarrow S_{23},~~~~~S_{34}\leftrightarrow 
S_{34},~~~~~S_{12}\leftrightarrow S_{21}, ~~~~~\mbox{etc.},\nonumber\\  
&\G^{1}_{~4}\leftrightarrow\G^{2}_{~4},~~~~~\G^{2}_{~2}
\leftrightarrow\G^{1}_{~1},~~~~~\mbox{etc.},\nonumber\\ 
&\Lambda\leftrightarrow\Lambda,~~~F\leftrightarrow F',~~~E
\leftrightarrow E',\\
&\Psi_\mu\leftrightarrow\Psi'_\mu, ~\mu=0,1,2,3,4,~~~~~~\phi_a\leftrightarrow
\phi'_a,~~a=0,1,2.
\eeq
It will be important that $\Lambda$ is invariant under the transformations 
(\ref{eq:trpr}).\\

\noindent
{\it 6) Distinguished totally null planes of maximal dimension.}\\

\noindent
Given a point $p\in\pp$ consider its corresponding pair $(s,s')$ of maximal
totally null planes of opposite selfduality in $\Pi (p)\in\cal M$. Lift $s$
horizontally to $\tilde{\tilde{s}}$ in $p$. It follows that
$\tilde{\tilde{s}}$ is a totally null 2-plane in any metric
$\tilde{\tilde{g}}$ in $\pp$.
There are only four 4-dimensional planes at $p$ which are totally null with 
respect to {\it any} 
metric $\tilde{\tilde{g}}$ and which 
contain $\tilde{\tilde{s}}$  as a subspace 
\footnote{Requirement that these planes must be null in any of metrics 
$\tilde{\tilde{g}}$ is crucial to obtain a discrete number of them!}. These 
may be defined  
as vector spaces $n_{F'EE'}$, $n_{F'E\bar{E}'}$, $n_{F'\bar{E}E'}$ and 
$n_{F'\bar{E}\bar{E}'}$ 
annihilating respectively $(\Lambda, F', E, E')$,
$(\Lambda, F', E, \bar{E}')$,
$(\Lambda, F', \bar{E}, E')$ and
$(\Lambda, F', \bar{E}, \bar{E}')$.
Point by point they define four maximal totally null bundles 
${\cal N}_{F'EE'}$, ${\cal N}_{F'E\bar{E}'}$, ${\cal N}_{F'\bar{E}E'}$ and 
${\cal N}_{F'\bar{E}\bar{E}'}$ over $\pp$. \\ 
Similarly, considering extensions of $\tilde{\tilde{s}}'$ we find four
other maximal totally null bundles
${\cal N}_{FEE'}$, ${\cal N}_{FE\bar{E}'}$, ${\cal N}_{F\bar{E}E'}$ and 
${\cal N}_{F\bar{E}\bar{E}'}$. Thus, in real cases, we will have eight 
different distinguished $r$-optical structures on $\pp$\\ 

\noindent
{\it 7) Invariant equations.}\\

\noindent
One set of geometrical equations on $\pp$ was already mentioned at the 
end of Section 5.3.4). By using forms $E$, $E'$, $\Lambda$, $F$ and $F'$ we 
can write further equations and systems of equations. Only those which 
are invariant under the 
allowed transformations of the forms have geometrical meaning. 
Among them there are the following.
\beq
&\der \Lambda\dz E\dz\Lambda\dz F=0\nonumber\\
&\der F\dz E\dz\Lambda\dz F=0\nonumber\\
&\der E \dz E\dz\Lambda\dz F=0\nonumber\\
&\der \Lambda\dz E'\dz\Lambda\dz F'=0\nonumber\\
&\der F'\dz E'\dz\Lambda\dz F'=0\nonumber\\
&\der E' \dz E'\dz\Lambda\dz F'=0\nonumber\\
&\der E\dz E\dz\Lambda\dz F'=0\nonumber\\
&\der E'\dz E'\dz\Lambda\dz F=0\nonumber\\
&\der E\dz E'\dz\Lambda\dz F'=0\nonumber\\
&\der E\dz E'\dz\Lambda\dz F=0.\nonumber
\eeq
One can continue this list. We discuss some of these equations 
in Section 6.\\

\section{Selfduality and the Einstein-Weyl equations}
Let $({\cal M}, g, A)$ be a 4-dimensional Weyl geometry. In this section we 
study relations between curvature properties of $({\cal M}, g, A)$ and 
the integrability conditions for the natural objects on the corresponding 
twistor bundles 
$\Pen$, $\Pen'$ and $\pp$. We start our analysis by giving a geometrical 
interpretation of the invariant equations on $\pp$ (see Section 4.3.7).
%
We use lthe ocal representation (\ref{eq:fonpp}) of the natural forms on $\pp$ 
to prove  
the following Lemma. 
\noindent
\bl~\\
The equations 
\be
{\rm d}\Lambda\dz\Lambda\dz F\dz E=0
\label{eq:integr1}
\ee
\be
{\rm d}F\dz\Lambda\dz F\dz E =0
\label{eq:integr2}
\ee
are automatically satisfied everywhere on $\Pen\Pen'$.
\el
Proof.\\
By using definitions (\ref{eq:fonpp}) and (\ref{eq:skr})
we easily compute that 
\beq
&\der\Lambda\dz \Lambda\dz F\dz E=\label{eq:wc1}\\
&[\Gamma^3_{~4}+z(\Gamma^2_{~4}-\Gamma^3_{~1})+
z'(\Gamma^1_{~4}-\Gamma^3_{~2})+zz'(\Gamma^3_{~3}+
\Gamma^4_{~4}-\Gamma^2_{~2}-\Gamma^1_{~1})-z^2\Gamma^2_{~1}-
z'^2\Gamma^1_{~2}+\nonumber\\
&zz'^2(\Gamma^1_{~3}-\Gamma^4_{~2})+
z'z^2(\Gamma^2_{~3}-\Gamma^4_{~1})+z^2z'^2\Gamma^4_{~3}]K\dz\Lambda
\dz F\dz E +\nonumber\\
&[z(\Gamma^2_{~2}+\Gamma^1_{~1}-\Gamma^3_{~3}-\Gamma^4_{~4})+
z'\Gamma^1_{~2}+zz'
(\Gamma^4_{~2}-\Gamma^1_{~3})-z^2z'\Gamma^4_{~3}]
F'\dz\Lambda\dz F\dz E\nonumber
\eeq
and
\beq
&\der F\dz \Lambda\dz F\dz E=
[z^2\Gamma^4_{~3}+z(\Gamma^1_{~3}-\Gamma^4_{~2})-
\Gamma^1_{~2}]F'\dz\Lambda\dz F\dz E+\nonumber\\ 
&[\Gamma^3_{~2}-\Gamma^1_{~4}+z'\Gamma^1_{~2}+z^2
(\Gamma^4_{~1}-\Gamma^2_{~3})+zz'(\Gamma^4_{~2}-\Gamma^1_{~3})-
z'z^2\Gamma^4_{~3})]K\dz\Lambda\dz F\dz E \label{eq:wc2}.
\eeq
Now, the right hand sides of the above expressions are actually 
equal to zero due to the Weyl geometry relations (\ref{eq:weylg}), which in 
the tetrad (\ref{eq:nt}) read 
$$\Gamma^1_{~2}=\Gamma^2_{~1}=\Gamma^3_{~4}=\Gamma^4_{~3}=0,$$
$$\Gamma^1_{~3}=\Gamma^4_{~2},~~~\Gamma^2_{~3}=\Gamma^4_{~1},~~~
\Gamma^1_{~4}=
\Gamma^3_{~2},~~~\Gamma^2_{~4}=\Gamma^3_{~1},$$
$$\Gamma^1_{~1}+\Gamma^2_{~2}=\Gamma^3_{~3}+\Gamma^4_{~4}=A.$$
This concludes the proof of the Lemma\footnote{
Note that the Weyl connections (\ref{eq:wlc1})-(\ref{eq:wlc}) are not 
the only connections that imply equations 
(\ref{eq:integr1})-(\ref{eq:integr2}). For our purposes, however, it is 
enough to restrict to Weyl connections.}. 

\noindent
Since equations (\ref{eq:integr1})-(\ref{eq:integr2}) are satisfied on 
$\Pen\Pen'$ it is natural to ask when the forms $\Lambda$, $F$ and $E$ form a 
closed differential ideal. Given the equations 
(\ref{eq:integr1})-(\ref{eq:integr2}), 
this question is equivalent to a question 
as to when an object $\der E\dz \Lambda\dz F\dz E$ vanishes identically 
on $\Pen\Pen'$. Also the related question of the vanishing of 
$\der E\dz\Lambda\dz F'\dz E$ may be interesting. A long but straightforward 
calculation leads to the following expressions
\be
\der E\dz\Lambda\dz F\dz E=[-\Psi_0+4\Psi_1z-6\Psi_2z^2
+4\Psi_3z^3-\Psi_4z^4]F'\dz K\dz\Lambda\dz F\dz E\label{eq:www}
\ee
\beq
&\der E\dz\Lambda\dz F'\dz E=
[-\frac{1}{2}S_{44} + z'S_{24} + 
zS_{14} - zz' (S_{12}+S_{34}) -\frac{1}{2} z'^2 S_{22}
\nonumber\\
&-\frac{1}{2} z^2 S_{11} + 
z'^2 zS_{23} + z^2z'S_{13} - 
\frac{1}{2}z^2z'^2S_{33}]F\dz K\dz\Lambda\dz F'\dz E.\label{eq:rrr} 
\eeq
Here we have used notations of Section 3 applied to null tetrad (\ref{eq:nt}).
The above formulae are implied by the general expressions for 
differentials of $\Lambda$, $F$, $F'$, $E$, $E'$ and $K$, which can be found 
in Appendix C.
\bt~\\
\vspace{-.8truecm}
\begin{itemize}
\item[1).] The equation
\be
{\rm d}E\dz\Lambda\dz F\dz E=0
\label{eq:weyl}
\ee
is satisfied everywhere on $\Pen\Pen'$ if and only if the underlying 
Weyl geometry $(g, A)$ on $\M$ is anti-self-dual.
\item[2).] The equation 
\be
{\rm d}E \dz\Lambda\dz F'\dz E =0
\label{eq:ric}
\ee
is satisfied everywhere on $\Pen\Pen'$ if and only if the underlying 
Weyl geometry $(g, A)$ satisfies the Einstein-Weyl equations 
on $\M$.
\end{itemize}
\et
Proof.\\
In the null tetrad (\ref{eq:nt})-(\ref{eq:met1})
condition (\ref{eq:slad}) expressing the zero trace of $S_{ij}$ reads 
$S_{12}-S_{34}=0$. Then the theorem follows immediately from 
expressions (\ref{eq:www})-(\ref{eq:rrr}) and the requirement that their 
right hand sides vanish for any $z$ and $z'$. $\Box$\\

\noindent
Applying the transformations (\ref{eq:trpr}) we also have analogous facts for 
primed objects.  

\bl~\\
The equations 
\be
{\rm d}\Lambda\dz\Lambda\dz F'\dz E'=0
\label{eq:integr1'}
\ee
\be
{\rm d}F'\dz\Lambda\dz F'\dz E' =0
\label{eq:integr2'}
\ee
are automatically satisfied everywhere on $\Pen\Pen'$.
\el
\noindent
\bt~\\
\vspace{-.8truecm}
\begin{itemize}
\item[1).] The equation
\be
{\rm d}E'\dz\Lambda\dz F'\dz E'=0
\label{eq:weyl'}
\ee
is satisfied everywhere on $\Pen\Pen'$ if and only if the underlying 
Weyl geometry $(g, A)$ on $\M$ is self-dual.
\item[2).] The equation 
\be
{\rm d}E' \dz\Lambda\dz F\dz E' =0
\label{eq:ric'}
\ee
is satisfied everywhere on $\Pen\Pen'$ if and only if the underlying 
Weyl geometry $(g, A)$ satisfies the Einstein-Weyl equations 
on $\M$.
\end{itemize}
\et
Transformations (\ref{eq:trpr}) show ($\Lambda$ is invariant!) 
that $dE\dz \Lambda\dz F'\dz E=0$ if and only if 
$dE'\dz \Lambda\dz F\dz E'=0$\footnote{Both equations here, are equivalent to 
the Einstein-Weyl 
equations for $({\cal M},g,A)$. This is due to the fact that the 
symmetric Ricci tensor $R_{(ij)}$ (or, equivalently $S_{ij}$ and $R$) 
is fully encoded in the differential of the spin connection $E$. The same 
information about the symmetric Ricci tensor is also encoded in d$E'$. To 
see this it is enough to note that $R$ as well as the quantity   
$$[-\frac{1}{2}S_{44} + z'S_{24} + 
zS_{14} - zz' (S_{12}+S_{34}) -\frac{1}{2} z'^2 S_{22}-
\frac{1}{2} z^2 S_{11} + 
z'^2 zS_{23} + z^2z'S_{13} - 
\frac{1}{2}z^2z'^2S_{33}]$$ 
of equation (\ref{eq:rrr}) are invariant under the transformations 
(\ref{eq:trpr}). See also Appendix C for explicit forms of d$E$ and d$E'$.}. 
This observation, together with the above theorems, leads 
to the following, interesting corollary.
\bc~\\
A Weyl geometry $({\cal M}, g, A)$ is anti-self-dual and satisfies the  
Einstein-Weyl equations if and only if the forms 
$(\Lambda, F, E, E')$ form a closed differential ideal  on $\pp$, 
i.e. iff on $\pp$ we have 
\beq
&\der \Lambda\dz E\dz\Lambda\dz F=0\nonumber\\
&\der F\dz E\dz\Lambda\dz F=0\nonumber\\
&\der E \dz E\dz\Lambda\dz F=0\nonumber\\
&\der E'\dz E'\dz\Lambda\dz F=0.\nonumber
\eeq
$({\cal M}, g, A)$ is self-dual and satisfies the  
Einstein-Weyl equations if and only if on $\pp$ we have 
\beq
&\der \Lambda\dz E'\dz\Lambda\dz F'=0\nonumber\\
&\der F'\dz E'\dz\Lambda\dz F'=0\nonumber\\
&\der E' \dz E'\dz\Lambda\dz F'=0\nonumber\\
&\der E\dz E\dz\Lambda\dz F'=0,\nonumber
\eeq
i.e. iff $(\Lambda, F', E, E')$ form a closed differential ideal  on $\pp$. 
\ec
The second part of the corollary follows from the first by applying 
(\ref{eq:trpr}).\\
According to Section 4.3.6 the forms $(\La , F, E, E')$ and 
$(\La, F', E, E')$ define natural maximal totally null bundles 
${\cal N}_{FEE'}$ 
and ${\cal N}_{F'EE'}$ on $\pp$. Using Lemma 3.1 we find that the 
above corollary has the following geometrical interpretation. 
\bt~\\
\vspace{-.8truecm}
\begin{itemize}
\item[1.)] The natural totally null bundle
${\cal N}_{FEE'}$ of maximal dimension 
over $\pp$ satisfies the integrability conditions (\ref{eq:integr}) 
if and only if the corresponding Weyl geometry $({\cal M}, g, A)$ is 
anti-self-dual and satisfies the Einstein-Weyl 
equations. 
\item[2.)] The natural totally null bundle  
${\cal N}_{F'EE'}$ of maximal dimension 
on $\pp$ satisfies the integrability conditions (\ref{eq:integr}) 
if and only if the corresponding Weyl geometry $({\cal M}, g, A)$ is 
self-dual and satisfies the Einstein-Weyl 
equations.
\end{itemize}
\et
Looking at the formula (\ref{eq:larz}) (see Appendix C), which gives a
differential $\der \La$, we see that
$\der\La\dz\La\dz F\dz\bar{E}\dz E'
=F'\dz E\dz\La\dz F\dz\bar{E}\dz E'\neq 0$. This means that the system  
$(\La, F, \bar{E},E')$ never forms a closed
differential ideal. Thus the totally null
bundle ${\cal N}_{F\bar{E}E'}$ has no chance to be integrable.
Similarly, since   
$\der\La\dz\La\dz F\dz\bar{E}\dz \bar{E}'
=F'\dz E\dz\La\dz F\dz\bar{E}\dz \bar{E}'\neq 0$, 
$\der\La\dz\La\dz F'\dz E\dz \bar{E}'
=F\dz E'\dz\La\dz F'\dz E\dz \bar{E}'\neq 0$, 
$\der\La\dz\La\dz F'\dz\bar{E}\dz \bar{E}'
=F\dz E'\dz\La\dz F'\dz\bar{E}\dz \bar{E}'\neq 0$ then also 
neither of the bundles 
${\cal N}_{F\bar{E}\bar{E}'}$,
${\cal N}_{F'E\bar{E}'}$, 
${\cal N}_{F'\bar{E}\bar{E}'}$ is integrable in the sense of definition 
(\ref{eq:integr}). To study the integrability conditions (\ref{eq:integr}) 
for bundles ${\cal N}_{FE\bar{E}'}$ and 
${\cal N}_{F'\bar{E}E'}$ we need to specify the real structure 
on $\cal M$. We postpone discussion of this case to the next Section.\\

\noindent
Analogous results about 
geometrical objects on $\Pen$ and $\Pen'$ (see Appendix D 
for differentials of 
the basis 1-forms) are summarized below.
\bt~\\
A Weyl geometry $({\cal M},g,A)$ is anti-self-dual if 
and only if the forms ($F, B, E$) form a closed differential ideal on $\Pen$,  
i.e. iff
\beq
&\der F\dz E\dz F\dz B=0\nonumber\\
&\der B\dz E\dz F\dz B=0\label{eq:ahps}\\
&\der E \dz E\dz F\dz B=0.\nonumber
\eeq
Similarly, a Weyl geometry $({\cal M},g,A)$ is self-dual iff on $\Pen'$ we 
have
\beq
&\der F'\dz E'\dz F'\dz B'=0\nonumber\\
&\der B'\dz E'\dz F'\dz B'=0\label{eq:ahps'}\\
&\der E'\dz E'\dz F'\dz B'=0.\nonumber
\eeq
\et
An obvious reinterpretation of this theorem in terms of the natural maximal 
totally null bundles on $\Pen$ and $\Pen'$ reads as follows. 
\bt~\\
A Weyl geometry $({\cal M},g,A)$ is anti-self-dual if 
and only if the natural totally null bundle ${\cal N}_E$ on $\Pen$ is  
integrable.\\
Similarly, a Weyl geometry $({\cal M},g,A)$ is self-dual iff on $\Pen'$ the 
natural totally null bundle ${\cal N}_{E'}$ is integrable.
\et
Looking at the differentials of $F$ on $\Pen$ and $F'$ on $\Pen'$
(see Appendix D) we find that
$\der F\dz F\dz B\dz \bar{E}=E\dz K\dz F\dz B\dz \bar{E}\neq 0$ and 
$\der F'\dz F'\dz B'\dz \bar{E'}=E'\dz K\dz F'\dz B'\dz \bar{E}'\neq 0$, 
respectively on $\Pen $ and $\Pen'$. This proves the following statement.
\bt~\\
Neither of the natural totally null bundles ${\cal N}_{\bar{E}}$ on $\Pen$ and 
${\cal N}_{\bar{E}'}$ on $\Pen'$ is integrable.
\et
The above two theorems are the Weyl-geometric counterparts of the 
Atiyah-Hitchin-Penrose-Singer 
theorems \cite{bi:AHS,bi:Pen3} for Lorentzian and Euclidean Riemannian 
4-manifolds. It is interesting that the integrability conditions of 
${\cal N}_E$ and ${\cal N}_{E'}$ say nothing about $A$. They only 
restrict the possible metrics on $\cal M$. 

\section{Real structures}
\noindent
{\bf 7.1 Reality conditions for the natural structures on twistor bundles}\\

\noindent  
In this section we consider real Weyl geometries $({\cal M},g,A)$. This  
means that 
the metric $g$ and the 1-form $A$ are real-valued. Such Weyl 
geometries and their twistor bundles 
are particular cases of the Weyl geometries considered in 
previous sections. Hence, our results of the previous sections are 
also valid here. In particular, a 
null cotetrad $(M,P,N,K)$ for $g$ may be chosen in such a way that 
\be
\bar{M}=P,~~~~\bar{N}=(1-|\varepsilon |)N+\varepsilon K,~~~~\bar{K}=
\varepsilon N +(1-|\varepsilon |) K,\label{eq:reco}
\ee
where $\varepsilon=0,1,-1$ for Lorentzian, neutral and Euclidean signature, 
respectively. Equations (\ref{eq:reco})  
imply the following reality conditions for the Weyl connection 1-forms.
\beq
&\bar{\G}^2_{~2}=\G^1_{~1}\nonumber\\
&\bar{\G}^2_{~3}=(1-|\varepsilon|)\G^1_{~3}+\varepsilon\G^1_{~4}\nonumber\\
&\bar{\G}^2_{~4}=(1-|\varepsilon|)\G^1_{~4}+\varepsilon\G^1_{~3}\nonumber\\
&\bar{\G}^1_{~3}=(1-|\varepsilon|)\G^2_{~3}+\varepsilon\G^2_{~4}
\label{eq:recocon}\\
&\bar{\G}^1_{~4}=(1-|\varepsilon|)\G^2_{~4}+\varepsilon\G^2_{~3}\nonumber\\
&\bar{\G}^4_{~4}=(1-2|\varepsilon|)\G^4_{~4}+|\varepsilon| A\nonumber\\
&\bar{\G}^3_{~3}=(1-2|\varepsilon|)\G^3_{~3}+|\varepsilon| A\nonumber
\eeq
Reality conditions for the curvature coefficients can be obtained from 
these equations.\\

\noindent
The above reality conditions imply the following properties of the natural 
maximal totally null bundles on $\Pen$, $\Pen'$ and $\pp$.\\

\noindent
{\it 1). Reality conditions for the maximal totally null bundles 
on $\Pen$ and $\Pen'$.}\\

\noindent
On $\Pen$ we have two natural maximal totally null bundles 
${\cal N}_{E}$ and ${\cal N}_{\bar{E}}$. 
It is easy to see that their real indices 
$r_E$ and $r_{\bar{E}}$ are equal. They  
depend on the signature of $g$ according to 
$r_E=r_{\bar{E}}=1-|\varepsilon|$. In the local representation 
$(x,z,\bar{z})$ of 
$\Pen$ one finds that
\be
{\cal N}_E\cap\bar{\cal N}_E={\cal N}_{\bar{E}}\cap\bar{\cal N}_{\bar{E}}=
(1-|\varepsilon|)
(\tilde{k}-z\tilde{m}-\bar{z}\tilde{p}+z\bar{z}\tilde{n}),
\ee
where $(\tilde{m},\tilde{p},\tilde{n},\tilde{k})$ are 
the horizontal lifts of the null tetrad $(m,p,n,k)$ from $\M$ to 
$\Pen$. An analogous formula for 
${\cal N}_{E'}$ and ${\cal N}_{\bar{E}'}$ reads  
\be
{\cal N}_{E'}\cap\bar{\cal N}_{E'}=
{\cal N}_{\bar{E}'}\cap\bar{\cal N}_{\bar{E}'}
=(1-|\varepsilon|)
(\tilde{k}-z'\tilde{p}-\bar{z}'\tilde{m}+z'\bar{z}'\tilde{n}).
\ee
{\it 2). Reality conditions for the maximal totally null bundles on $\pp$.}\\
We have eight natural totally null bundles on $\pp$. Their reality conditions 
are given below in the local representation $(x,z,\bar{z},z',\bar{z}')$ of 
$\pp$.
\beq
&{\cal N}_{FEE'}\cap\bar{\cal N}_{FEE'}={\cal N}_{FE\bar{E}'}
\cap\bar{\cal N}_{FE\bar{E}'}=
{\cal N}_{F\bar{E}E'}\cap\bar{\cal N}_{F\bar{E}E'}=
{\cal N}_{F\bar{E}\bar{E}'}
\cap\bar{\cal N}_{F\bar{E}\bar{E}'}=\nonumber\\
&=(1-|\varepsilon|)
(\tilde{\tilde{k}}-z\tilde{\tilde{m}}-
\bar{z}\tilde{\tilde{p}}+z\bar{z}\tilde{\tilde{n}}),\nonumber\\
&{\cal N}_{F'EE'}\cap\bar{\cal N}_{F'EE'}={\cal N}_{F'E\bar{E}'}
\cap\bar{\cal N}_{F'E\bar{E}'}=
{\cal N}_{F'\bar{E}E'}\cap\bar{\cal N}_{F'\bar{E}E'}=
{\cal N}_{F'\bar{E}\bar{E}'}
\cap\bar{\cal N}_{F'\bar{E}\bar{E}'}=\nonumber\\
&=(1-|\varepsilon|)
(\tilde{\tilde{k}}-z'\tilde{\tilde{p}}-
\bar{z}'\tilde{\tilde{m}}+z'\bar{z}'\tilde{\tilde{n}}).\nonumber
\eeq
Here $(\tilde{\tilde{m}},\tilde{\tilde{p}},\tilde{\tilde{n}},
\tilde{\tilde{k}})$ are 
horizontal lifts of null tetrad $(m,p,n,k)$ to 
$\pp$. \\

\noindent
Note, that the real indices of all of the natural maximal totally null bundles
on $\Pen$, $\Pen'$ and $\pp$ are either 0 or 1. Thus, the
twistor bundles get naturally equipped with either
Hermitian or optical geometries. We do not see possibilities of
distinguishing $r$-optical geometries with $r>1$ on
$\Pen$, $\Pen'$ and $\pp$.\\

\noindent
{\bf 7.2 Euclidean and neutral signature}\\

\noindent
These cases are characterized by $|\varepsilon|=1$. It follows from 
Section 7.1 that the real indices of 
all the natural maximal totally null bundles are equal to zero. This, in 
particular, means   
that ${\cal N}_E$ and ${\cal N}_{\bar{E}}$ define almost complex 
structures ${\cal J}_E$ and  ${\cal J}_{\bar{E}}$ on $\Pen$.
These structures are almost hermitian
in any metric $\tilde{g}$. Similarly, we have two natural almost hermitian   
structures $({\cal J}_{E'}, \tilde{g}')$ and 
$({\cal J}_{\bar{E}'},\tilde{g}')$ on $\Pen'$, and eight almost hermitian  
structures 
$({\cal J}_{FEE'},\tilde{\tilde{g}})$, 
$({\cal J}_{FE\bar{E}'},\tilde{\tilde{g}})$, 
$({\cal J}_{F\bar{E}E'},\tilde{\tilde{g}})$, 
$({\cal J}_{F\bar{E}\bar{E}'},\tilde{\tilde{g}})$, 
$({\cal J}_{F'EE'},\tilde{\tilde{g}})$, 
$({\cal J}_{F'E\bar{E}'},\tilde{\tilde{g}})$, 
$({\cal J}_{F'\bar{E}E'},\tilde{\tilde{g}})$,   
$({\cal J}_{F'\bar{E}\bar{E}'},\tilde{\tilde{g}})$ on $\pp$.  
Integrability conditions of these almost hermitian structures are  
equivalent to conditions (\ref{eq:integr}) for the corresponding 
${\cal N}$s. Most of them were already studied in Section 6.
The integrability conditions for the remaining two
structures $({\cal J}_{FE\bar{E}'},\tilde{\tilde{g}})$ and
$({\cal J}_{F'\bar{E}E'},\tilde{\tilde{g}})$
follow from the expressions of Appendix C. In particular, looking at
(\ref{eq:larz}) and the primed counterpart of (\ref{eq:frz}) we find that
$$\der\La\dz\La\dz F'\dz E'\dz\bar{E}\equiv 0,$$
$$\der F'\dz\La\dz F'\dz E'\dz\bar{E}\equiv 0.$$
On the other hand (\ref{eq:prcp}) and the primed counterpart of
(\ref{eq:dere}) show that
$$\der E'\dz\La\dz F'\dz E'\dz\bar{E}\equiv 0,$$
$$\der\bar{E}\dz\La\dz F'\dz E'\dz\bar{E}\equiv 0$$
if and only if the Weyl geometry is self-dual and satisfies the
Einstein-Weyl equations. Thus we find that
$({\cal J}_{F'\bar{E}E'},\tilde{\tilde{g}})$ is integrable
only for such Weyl geometries.
An analogous result holds also for
$({\cal J}_{FE\bar{E}'},\tilde{\tilde{g}})$. This leads to the
following theorem.
\bt~\\
Let $\cal M$ be a 4-dimensional real manifold equipped with a Euclidean or  
a neutral signature 
Weyl geometry $(g, A)$. Let $\Pen$ (respectively, $\Pen'$) be a 
corresponding twistor bundle of all self-dual 
(respectively, anti-self-dual) maximal totally null spaces over 
$\cal M$. Let $\pp$ be a fiber product of the bundles $\Pen$ and $\Pen'$.
\begin{itemize}
\item[1.)] 
There are two 
natural almost hermitian structures ${\cal J}_{E}$ and 
${\cal J}_{\bar{E}}$ on $\Pen$. ${\cal J}_E$ is integrable if and only if 
the Weyl geometry is anti-self-dual. 
${\cal J}_{\bar{E}}$ is never integrable.
\item[2.)]
There are two 
natural almost hermitian structures ${\cal J}_{E'}$ and 
${\cal J}_{\bar{E}'}$ on $\Pen'$. ${\cal J}_E'$ is integrable if and only 
if the Weyl geometry is self-dual. ${\cal J}_{\bar{E}'}$ is never integrable.
\item[3.)]
There ~are ~eight ~natural ~almost hermitian structures 
$({\cal J}_{FEE'},\tilde{\tilde{g}})$, 
$({\cal J}_{FE\bar{E}'},\tilde{\tilde{g}})$, 
$({\cal J}_{F\bar{E}E'},\tilde{\tilde{g}})$, 
$({\cal J}_{F\bar{E}\bar{E}'},\tilde{\tilde{g}})$, 
$({\cal J}_{F'EE'},\tilde{\tilde{g}})$, 
$({\cal J}_{F'E\bar{E}'},\tilde{\tilde{g}})$, 
$({\cal J}_{F'\bar{E}E'},\tilde{\tilde{g}})$,   
$({\cal J}_{F'\bar{E}\bar{E}'},\tilde{\tilde{g}})$ on $\pp$.
\begin{itemize}
\item[i.)]
$({\cal J}_{FEE'},\tilde{\tilde{g}})$  
is integrable if and only if  
the Weyl geometry is anti-self-dual and satisfies the Einstein-Weyl 
equations. These integrability conditions are also equivalent to the 
integrability of $({\cal J}_{FE\bar{E}'},\tilde{\tilde{g}})$. 
\item[ii.)]
$({\cal J}_{F'EE'},\tilde{\tilde{g}})$  
is integrable if and only if  
the Weyl geometry is self-dual and satisfies the Einstein-Weyl 
equations. These integrability conditions are also equivalent to the 
integrability of $({\cal J}_{F\bar{E}E'},\tilde{\tilde{g}})$. 
\item[iii.)]
$({\cal J}_{F\bar{E}E'},\tilde{\tilde{g}})$, 
$({\cal J}_{F\bar{E}\bar{E}'},\tilde{\tilde{g}})$, 
$({\cal J}_{F'E\bar{E}'},\tilde{\tilde{g}})$, 
$({\cal J}_{F'\bar{E}\bar{E}'},\tilde{\tilde{g}})$ are never integrable.
\end{itemize}
\end{itemize}
\et
Note that there are no restrictions on the potential $A$ in this theorem. 
This is a generalization of the classical Atiyah-Hitchin-Singer theorem 
\cite{bi:AHS}. 
It is interesting that in a purely Euclidean case 
($\varepsilon=-1$, $A=0$) we have a 
holomorphic interpretation of the self-dual Einstein equations.\\

In Section 6, Theorem 6.1, we interpreted an invariant equation 
\be
\der E\dz \La\dz F'\dz E=0 \label{eq:ieew}
\ee
on $\pp$ as a necessary and sufficient condition for the 
Weyl geometry to satisfy the Einstein-Weyl equations. In the present real 
case we can reformulate this fact in the holomorphic language.\\ 
Given an almost complex structure $\cal J$ on $\pp$ we can decompose 
the complexification T$(\pp)^{\bf C}$ of its tangent bundle onto the 
eigenspaces of $\cal J$. The $+i$, $-i$ eigenspaces are denoted respectively 
by T$^{(1,0)}$ and T$^{(0,1)}$. One easily finds that T$^{(0,1)}$ is 
the same as the maximal totally null  
bundle $\cal N$ representing $\cal J$. The above decomposition of 
T$(\pp)^{\bf C}$ induces an analogous decomposition of 
the complexification T$^*(\pp)^{\bf C}$ of the cotangent space.. We denote 
by T$^{*(1,0)}$ the anihilator of T$^{(0,1)}$ and by 
T$^{*(0,1)}$ the anihilator of T$^{(1,0)}$. In an analogous way 
T$^{*(u,w)}$ denotes an exterior product of $u$ copies of T$^{*(1,0)}$  
and $w$ copies of T$^{*(0,1)}$ bundles. It is well known that we have the 
following decomposition of the bundle $\La^2$ of all complex-valued 2-forms
on
$\pp$
\be
\La^2={\rm T}^{*(2,0)}\oplus{\rm T}^{*(1,1)}\oplus{\rm T}^
{*(0,2)}.\label{eq:rozklad}
\ee
Consider now a natural almost complex structure ${\cal J}_{F'EE'}$ on $\pp$ 
and its corresponding decomposition of $\La^2$. We analyze the differential 
$\der E$ of the spin connection 1-form from the point of view of this 
decomposition. To do this we consider a set $\cal W$ of 2-forms over
$\pp$ defined by
$${\cal W}=\{ w\in\G(\La^2) ~~{\rm s.t.}~~
w=~\Omega\dz E+t~\der E\},$$
where $\Omega$ and $t$ are respectively any complex-valued 1-form and
function on $\pp$.
We decompose $\cal W$ according to (\ref{eq:rozklad}).
\bt~\\
A real Euclidean- or a neutral-signature  4-dimensional 
Weyl geometry $({\cal M}, g, A)$ satisfies the Einstein-Weyl 
equations if and only if 
$${\cal W}\cap\G({\rm T}^{*(0,2)})=\{0\}.$$
\et
Proof.\\
Observe that sections of 
T$^{*(2,0)}$ and T$^{*(1,1)}$ are in the ideal generated by the forms 
$(\La, F', E, E')$. On the other hand the set of sections of
T$^{*(0,2)}$ has zero
intersection with this ideal. Thus in the decomposition of 
the sections of T$^{*(0,2)}$ onto a basis corresponding to 
$(F, F', \La, T, E, \bar{E}, E', \bar{E}')$ there are no 
forms $\La$, $F'$, $E$, $E'$. 
Now the proof follows directly from the differential of $E$ given
by (\ref{eq:dere}) (see Appendix C).\\

We close this section by considering a Weyl geometry which is not
anti-self-dual. Fibers of its twistor bundle
$\Pen$ have a discrete number of distinguished points. To see this
consider such points in a fiber $\pi^{-1}(x)$, 
$x\in\cal M$, in which the expression $\der E\dz E\dz F\dz B$ vanishes. Due to
\be
\der E\dz E\dz F\dz B
=[-\Psi_0+4\Psi_1z-6\Psi_2z^2
+4\Psi_3z^3-\Psi_4z^4]P\dz K\dz E\dz F\dz B\label{eq:81}
\ee
we find that in the not anti-self-dual case there are at most 
four $z$s, corresponding to four points at the fiber $\pi^{-1}(x)$, 
in which the right hand side of (\ref{eq:81}) is zero.
These four points correspond to four maximal
totally null self-dual planes at $x$. Thus, in a generic case, at every 
point of $\cal M$ we have four distinguished almost hermitian structures. It 
further follows from reality conditions for $\Psi_\mu$ that in the not 
anti-self-dual case these four structures group in pairs of 
mutually conjugated structures. These two pairs may 
coincide in particular cases and, together with the additional 
two pairs associated  
with the similar considerations on $\Pen'$, may be used to  
classify Weyl geometries. An interesting fact is that in the 
non-half-flat case these distinguished almost hermitian structures are the 
only ones that may be integrable on $\cal M$ \cite{bi:Derdz}. A less well 
known 
fact is that in a purely Euclidean case ($\varepsilon=-1, A=0)$, if the  
Einstein equations are satisfied, then any of the 
distinguished almost hermitian 
structures is integrable \cite{bi:GS,bi:phd,bi:Przan}.\\

\noindent
{\bf 7.3 Lorentzian signature}\\

\noindent
Due to the condition
$\varepsilon=0$ 
the fibers of  
${\cal N}_{FEE'}\cap\bar{\cal N}_{FEE'}$, ${\cal N}_{FE\bar{E}'}
\cap\bar{\cal N}_{FE\bar{E}'}$, 
${\cal N}_{F\bar{E}E'}\cap\bar{\cal N}_{F\bar{E}E'}$ and  
${\cal N}_{F\bar{E}\bar{E}'}
\cap\bar{\cal N}_{F\bar{E}\bar{E}'}$ are all 1-dimensional and, in every  
point $p\in\pp$, are spanned by a real vector 
$$\kappa=\tilde{\tilde{k}}-z\tilde{\tilde{m}}-
\bar{z}\tilde{\tilde{p}}+z\bar{z}\tilde{\tilde{n}}.$$
The fibers of  
${\cal N}_{F'EE'}\cap\bar{\cal N}_{F'EE'}$, ${\cal N}_{F'E\bar{E}'}
\cap\bar{\cal N}_{F'E\bar{E}'}$, 
${\cal N}_{F'\bar{E}E'}\cap\bar{\cal N}_{F'\bar{E}E'}$, 
${\cal N}_{F'\bar{E}\bar{E}'}
\cap\bar{\cal N}_{F'\bar{E}\bar{E}'}$ are spanned by 
a real vector 
$$\kappa'=(\tilde{\tilde{k}}-z'\tilde{\tilde{p}}-
\bar{z}'\tilde{\tilde{m}}+z'\bar{z}'\tilde{\tilde{n}}).$$
These two real vectors are null and, together with their 
corresponding maximal totally null spaces, define eight distinguished optical 
geometries 
${\cal O}_{FEE'}$, ${\cal O}_{F\bar{E}E'}$, 
${\cal O}_{FE\bar{E}'}$, ${\cal O}_{F\bar{E}\bar{E}'}$, 
${\cal O}_{F'EE'}$, ${\cal O}_{F'\bar{E}E'}$, 
${\cal O}_{F'E\bar{E}'}$, ${\cal O}_{F'\bar{E}\bar{E}'}$ on $\pp$. According 
to Section 2 any  
7-dimensional sub manifold of $\pp$ transversal to $\kappa$ is 
equipped with four CR-structures that correspond to 
${\cal O}_{FEE'}$, ${\cal O}_{F\bar{E}E'}$, 
${\cal O}_{FE\bar{E}'}$, ${\cal O}_{F\bar{E}\bar{E}'}$. 
Another four 
families  
of CR-structures are associated with  
the 7-dimensional submanifolds of $\pp$ transversal to $\kappa'$. To get a
Lorentzian version of Theorem 7.1 we still need to note that if $\eps=0$
then $\Psi'_\mu=\bar{\Psi}_\mu$ and $\phi'_a=\bar{\phi}_a$ for
all $\mu=0,1,2,3,4$ and all $a=0,1,2$. Thus in this case (anti)-self-duality
of the Weyl tensor is equivalent to conformal flatness of the metric and
(anti)-self-duality of $\cal F$ means that ${\cal F}\equiv 0$.
\bt~\\
Let $\cal M$ be a 4-dimensional real manifold equipped with a Lorentzian   
Weyl geometry $(g, A)$. Let $\Pen$ (respectively, $\Pen'$) be the
corresponding twistor bundle of all self-dual
(respectively, anti-self-dual) maximal totally null spaces over 
$\cal M$. Let $\pp$ be a fiber product of bundles $\Pen$ and $\Pen'$.
\begin{itemize}
\item[1.)]
There are eight natural almost optical geometries 
${\cal O}_{FEE'}$, ${\cal O}_{FE\bar{E}'}$, 
${\cal O}_{F\bar{E}E'}$, ${\cal O}_{F\bar{E}\bar{E}'}$, 
${\cal O}_{F'EE'}$, 
${\cal O}_{F'E\bar{E}'}$, ${\cal O}_{F'\bar{E}E'}$, 
${\cal O}_{F'\bar{E}\bar{E}'}$ on $\pp$.
\item[2.)]
The following conditions are equivalent.
\begin{itemize}
\item[i.)]
${\cal O}_{FEE'}$  (respectively, ${\cal O}_{F'EE'}$) is integrable 
\item[ii.)]
there exists a unique integrable 7-dimensional CR-structure 
obtained from $\pp$ by identifying points lying on the same integral curve 
of $\kappa$ 
(respectively, $\kappa'$) and 
associated with ${\cal O}_{FEE'}$ (respectively, ${\cal O}_{F'EE'}$) 
\item[iii.)]
the Weyl geometry is conformally flat and satisfies the Einstein-Weyl 
equations. 
\end{itemize}
\item[3)]
The following conditions are equivalent.
\begin{itemize}
\item[i.)]
${\cal O}_{FE\bar{E}'}$  (respectively, ${\cal O}_{F'\bar{E}E'}$) is 
integrable 
\item[ii.)]
there exists a unique integrable 7-dimensional CR-structure 
obtained from $\pp$ by identifying points lying on the same integral curve 
of $\kappa$ 
(respectively, $\kappa'$) and associated with ${\cal O}_{FE\bar{E}'}$ 
(respectively, ${\cal O}_{F'\bar{E}E'}$) 
\item[iii.)]
the Weyl geometry is conformally flat, has zero scalar curvature $R$ and
a potential $A$ such that $\der A\equiv 0$.
\end{itemize}
\item[4.)]
${\cal O}_{F\bar{E}E'}$, ${\cal O}_{F\bar{E}\bar{E}'}$,   
${\cal O}_{F'E\bar{E}'}$ and ${\cal O}_{F\bar{E}\bar{E}'}$ are never 
integrable.  
\end{itemize}
\et

\noindent
Only point ${\em 3 iii)}$ of the theorem requires justification.
It is, however, easy to see that this follows from the differentials of
the basis 1-forms given in Appendix C and
from formula (\ref{eq:prcp}) with $\eps=0$.\\

\noindent
To find a passage from the above results to the description on the original 
Penrose bundle given in Section 2 we do as follows.\\ 
First, we note 
that $\Pen$ defined as a bundle of all self-dual 
totally null 2-planes in 
$\cal M$ is naturally isomorphic to the Penrose bundle of 
all real null directions in $\cal M$. This is due to the fact that in the 
case of 
a Lorentzian metric any 
2-dimensional totally null plane in $\cal M$ is in one-to-one correspondence  
with a null direction.\\
Second, we note that the natural field of real null 
directions $\kappa$ 
can be naturally pushed forward from $\pp$ to $\Pen$ by   
means of projection pr$_*$. This defines a field of real directions 
$X_L=$pr$_*\kappa$ on 
$\Pen$ which is null in any metric $\tilde{g}$. Using $X_L$ we define 
a field of directions of 
a real 1-form $\La_L=\tilde{g}(X_L)$. Now, we can take $E_L=E$ and 
define 1-forms $F_L$ (complex) 
and $T$ (real) by $\tilde{g}=h\bar{h}E_L\bar{E}_L+\La_LT+F_L\bar{F}_L$. It is 
easy to see that the forms $(\La_L,F_L,E_L)$ are precisely given up to 
transformations (\ref{eq:trla})-(\ref{eq:trf}). The above information 
is sufficient to reconstruct on $\Pen$ all the structures of Section 2. In 
particular, the anihilator ${\cal N}_{FE}$ of $(\La_L,F_L,E_L)$ 
defines an almost optical geometry ${\cal O}_{FE}$ on $\Pen$. Similarly 
the annihilators 
${\cal N}_{\bar{F}E}$, ${\cal N}_{F\bar{E}}$ and  
${\cal N}_{\bar{F}\bar{E}}$ 
of respectively  
$(\La_L,\bar{F}_L,E_L)$, $(\La_L,F_L,\bar{E}_L)$ and 
$(\La_L,\bar{F}_L,\bar{E}_L)$ define optical geometries 
${\cal O}_{\bar{F}E}$, ${\cal O}_{F\bar{E}}$ and  
${\cal O}_{\bar{F}\bar{E}}$. Their integrability conditions are summarized 
in the following Theorem.  
\bt~\\
There are four natural optical geometries 
${\cal O}_{FE}$, ${\cal O}_{\bar{F}E}$, ${\cal O}_{F\bar{E}}$ and  
${\cal O}_{\bar{F}\bar{E}}$ on a 6-dimensional Penrose twistor bundle of 
all null directions over a 4-dimensional manifold ${\cal M}$ equipped 
with a Lorentzian Weyl geometry $(g, A)$. 
\begin{itemize}
\item[1.)] The following conditions are equivalent.
\begin{itemize}
\item[i.)] ${\cal O}_{FE}$ or ${\cal O}_{\bar{F}\bar{E}}$ is integrable 
\item[ii.)] a 5-dimensional manifold of lines of a congruence generated by 
$X_L$ is naturally equipped with an integrable CR-structure
\item[iii.)] the metric $g$ is conformally flat.
\end{itemize}
\item[2.)] ${\cal O}_{\bar{F}E}$ and ${\cal O}_{F\bar{E}}$ are never 
integrable.
\end{itemize}
\et
We failed in finding a geometrical interpretation of the following theorem.
\bt~\\
A Lorentzian 4-dimensional Weyl geometry satisfies the Einstein-Weyl 
equations if the forms $(E_L,\bar{F}_L,\La_L)$ satisfy an invariant 
equation  
$$
\der E\dz\bar{F}_L\dz\La_L\dz E_L=0
$$
everywhere on $\Pen$.
\et
\section{Acknowledgements}
I am very grateful to A. Trautman who directed my attention to the topics 
presented in this paper. His constant guidance during the preparation of the  
paper was crucial. Many of the results presented here arose during 
discussions with him. I also wish to thank S. Bazanski, T. Friedrich, 
J. Kijowski, P. Kobak, W. Kopczynski, J. Lewandowski, N. Manojlovic, 
L. Mason, J. Mourao, E. T. Newman, R. Penrose, D. C. Robinson, 
S. M. Salamon, G. A. J. Sparling, J. Tafel, P. Tod and 
H. Urbantke for helpful discussions.\\
This work was performed during my stay in SISSA, Humboldt 
University in Berlin, 
University of Algarve in Portugal and King's College in London. I wish to 
thank 
C. Reina, T. Friedrich, N. Manojlovic and D. C. Robinson for financial 
support during these visits.   
\section{Appendix A}
In this Appendix we present formulae, which give decompositions
of certain Weyl-geometric objects
onto the Levi-Civita and $A$-dependent parts.\\
For the Riemann tensor 2-forms $\Omega_{ij}$ we have
\be
\Omega_{ij}=\omega_{ij}+\frac{1}{2}g_{ij}{\cal F}+DA_{[j}
\dz\theta_{i]}+\frac{1}{2}A_{[j}\theta_{i]}\dz A-\frac{1}{4}A^2
\theta_i\dz\theta_j,
\ee  
where $\om_{ij}$ denote curvature 2-forms of the Levi-Civita
connection $\g_{ij}$,
\be
DA_i=A_{i;j}\theta^j=
\der A_i - A_j\gamma^j_{~i}~~~~~~~~~\theta_i=g_{ik}\theta^k
\ee
and $A^2=g^{ij}A_iA_j$. 
From these equations one easily gets curvature coefficients 
\be
R_{ijkl}=r_{ijkl}+\frac{1}{2}g_{ij}{\cal F}_{kl}+
(A_{j;[k}g_{l]i}-A_{i;[k}g_{l]j})+\frac{1}{2}(A_jg_{i[k}A_{l]}-
A_ig_{j[k}A_{l]})-\frac{1}{2}A^2g_{i[k}g_{l]j}.
\ee
Here $r_{ijkl}$ are the usual components of the curvature tensor for 
the Levi-Civita connection $\g^i_{~j}$.\\
The Ricci tensor and Ricci scalar decompositions read respectively
\be
R_{jl}=
r_{jl}+\frac{1}{2}{\cal F}_{jl}-A_{j;l}-\frac{1}{2}g^{ik}A_{i;k}g_{jl}+
\frac{1}{2}A_jA_l-\frac{1}{2}A^2g_{jl},
\ee
\be
R=r-3g^{jl}A_{j;l}-\frac{3}{2}A^2.
\ee
Here quantities $r_{ij}$ and $r$ denote respectively 
the Ricci tensor and Ricci scalar of the Levi-Civita 
connection associated with $g$.

\noindent 
Symmetric part of the Ricci tensor decomposes according to 
\be
R_{(jl)}=
r_{jl}-A_{(j;l)}-\frac{1}{2}g^{ik}A_{i;k}g_{jl}+
\frac{1}{2}A_jA_l-\frac{1}{2}A^2g_{jl}.
\ee 
$S_{ij}$ has the following decomposition into the Levi-Civita ($s_{ij})$ 
and the $A$-dependent ($\sigma_{ij}$) part.
\be
S_{ij}=s_{ij}+\sigma_{ij},
\ee
where
\be
s_{ij}=r_{ij}-\frac{1}{4}rg_{ij}
\ee
and
\be
\sigma_{ij}=-A_{(i;j)}+\frac{1}{4}g^{kl}A_{k;l}g_{ij}+\frac{1}{2}A_iA_j-
\frac{1}{8}A^2g_{ij}.
\ee
Decomposition of $C_{ijkl}$ is given by
\be
C_{ijkl}=w_{ijkl}+f_{ijkl},
\ee
where the Levi-Civita ($w_{ijkl}$) and $A$-dependent ($f_{ijkl}$) 
parts read respectively
\be
w_{ijkl}=r_{ijkl}+\frac{1}{3}rg_{i[k}g_{l]j}+r_{j[k}g_{l]i}+r_{i[l}g_{k]j},
\ee
\be
f_{ijkl}=\frac{1}{2}({\cal F}_{j[k}g_{l]i}+{\cal F}_{i[l}g_{k]j})+
\frac{1}{2}g_{ij}{\cal F}_{kl}.
\ee
Returning to the curvature forms $\Omega_{ij}$ we decompose it onto 
antisymmetric and symmetric parts. We 
note that $\Omega_{[ij]}$ can be further decomposed onto a part 
$\Omega^U_{[ij]}$ with coefficients having all the symmetries of 
the usual (i.e. Levi-Civita connection) Riemann tensor and the remaining 
part 
$\Omega^{NU}_{[ij]}$. Explicitely we have 
\be
\Omega_{ij}=\Omega^U_{[ij]}+\Omega^{NU}_{[ij]}+\Omega_{(ij)},
\ee
where
\be
\Omega_{(ij)}=\frac{1}{2}g_{ij}{\cal F},
\ee
\be
\Omega^{NU}_{[ij]}=-\frac{1}{4}(\theta_i\dz{\cal F}_j-\theta_j\dz{\cal F}_i),
\ee
\be
\Omega^U_{[ij]}=w_{ij}+\frac{1}{12}R\theta_i\dz\theta_j+\frac{1}{2}
(\theta_i\dz S_j-\theta_j\dz S_i),
\ee
and ${\cal F}_i=\theta^k{\cal F}_{ki}$, $S_i=\theta^k S_{ki}$.

\section{Appendix B}

Given a Weyl geometry $({\cal M},g,A)$ consider null cotetrad (\ref{eq:nt}).
Then the decomposition of the antisymmetric
part of the curvature 2-forms $\Omega_{[ij]}$ onto basis of self-dual
and anti-self-dual 2-forms read as follows.
\beq
&&\Omega_{[14]}=\nonumber\\
&&\Psi'_0M\dz K+ (\Psi'_1+\frac{1}{4}\phi'_2)(N\dz K+M\dz P)+(\Psi'_2+
\frac{R}{12}+\frac{1}{2}\phi'_1)N\dz P+\nonumber\\
&&\frac{1}{2}S_{44}P\dz K+
\frac{1}{2}S_{41}(N\dz K-M\dz P)+\frac{1}{2}S_{11}N\dz M\nonumber\\
&&\nonumber\\
&&\Omega_{[23]}=\nonumber\\
&&(-\Psi'_2+\frac{R}{12}+\frac{1}{2}\phi'_1)M\dz K+(-\Psi'_3+
\frac{1}{4}\phi'_0)(N\dz K+M\dz P)-\Psi'_4N\dz P+\nonumber\\
&&-\frac{1}{2}S_{22} P\dz K -\frac{1}{2}S_{32}(N\dz K-M\dz P)-
\frac{1}{2}S_{33}N\dz M\label{eq:om'}\\
&&\nonumber\\
&&\frac{1}{2}(\Omega_{[34]}+\Omega_{[12]})=\nonumber\\
&&(\Psi'_1-\frac{1}{4}\phi'_2)M\dz K+ (\Psi'_2-\frac{R}{24})(N\dz K+M\dz P)+
(\Psi'_3+\frac{1}{4}\phi'_0)N\dz P+\nonumber\\
&&\frac{1}{2}S_{42} P\dz K+ \frac{1}{4}(S_{12}+S_{34})
(N\dz K-M\dz P)+\frac{1}{2}S_{31}N\dz M+\nonumber
\eeq

\beq
&&\Omega_{[24]}=\nonumber\\
&&\Psi_0P\dz K+ (\Psi_1+\frac{1}{4}\phi_2)(N\dz K-M\dz P)+(\Psi_2+
\frac{R}{12}+\frac{1}{2}\phi_1)N\dz M+\nonumber\\
&&\frac{1}{2}S_{44}M\dz K+
\frac{1}{2}S_{42}(N\dz K+M\dz P)+\frac{1}{2}S_{22}N\dz P\nonumber\\
&&\nonumber\\
&&\Omega_{[13]}=\nonumber\\
&&(-\Psi_2+\frac{R}{12}+\frac{1}{2}\phi_1)P\dz K+(-\Psi_3+
\frac{1}{4}\phi_0)(N\dz K-M\dz P)-\Psi_4N\dz M+\nonumber\\
&&-\frac{1}{2}S_{11} M\dz K -\frac{1}{2}S_{31}(N\dz K+M\dz P)-
\frac{1}{2}S_{33}N\dz P\label{eq:om''}\\
&&\nonumber\\
&&\frac{1}{2}(\Omega_{[34]}-\Omega_{[12]})=\nonumber\\
&&(\Psi_1-\frac{1}{4}\phi_2)P\dz K+ (\Psi_2-\frac{R}{24})(N\dz K-M\dz P)+
(\Psi_3+\frac{1}{4}\phi_0)N\dz M+\nonumber\\
&&\frac{1}{2}S_{41} M\dz K+ \frac{1}{4}(S_{12}+S_{34})
(N\dz K+M\dz P)+\frac{1}{2}S_{32}N\dz P+\nonumber
\eeq
\section{Appendix C}
In this Appendix we give differentials of basis 1-forms
$(E,\bar{E},E',\bar{E}',F,F',\La,T)$ on $\pp$.
For a given Weyl geometry $({\cal M},g,A)$ we use a null cotetrad
(\ref{eq:nt}) and represent the basis 1-forms on $\pp$ according to
(\ref{eq:fonpp}).\\
Let 
$$
2\g=\G_{211}+\G_{341}+2z\G_{131}-z'(\G_{213}+\G_{343})-2zz'\G_{133},
$$
\beq
&2\om=\G_{214}+\G_{344}+z(2\G_{134}-\G_{211}-\G_{341})-
z'(\G_{212}+\G_{342})\nonumber\\
&+zz'(\G_{213}+\G_{343}-2\G_{132})-2z^2\G_{131}+2z^2z'\G_{133},\nonumber
\eeq
$$
2\Phi=\frac{1}{2}S_{44} - z'S_{24} -
zS_{14}  + zz' (S_{12}+S_{34}) +\frac{1}{2} z'^2 S_{22}
+\frac{1}{2} z^2 S_{11} - z'^2 zS_{23} - z^2z'S_{13} + 
\frac{1}{2}z^2z'^2S_{33}, 
$$
$$
4\Psi=\Psi_0-4\Psi_1z+6\Psi_2z^2-4\Psi_3z^3+\Psi_4z^4,
$$
$$
4\phi=-\phi_2+2z\phi_1-z^2\phi_0,
$$
$$
2a=A_1z+A_2z'-A_3zz'-A_4.
$$
Applying transformations (\ref{eq:trpr}) we get also quantities $\g'$,
$\om'$, $\Psi'$ and $\phi'$
\footnote{Note that these transformations do not affect $\Phi$ and $a$.}.\\
Vanishing of some of the above coefficients has well defined meaning in terms
of the Weyl geometry on $\cal M$. Note, for example, that $\Phi\equiv 0$
means that the Weyl geometry satisfies the Einstein-Weyl equations,
$\Psi\equiv 0$ means that the Weyl geometry is anti-self-dual,
$\phi'\equiv 0$ means that curvature $\cal F$ of $A$ is self-dual.\\
Using the above quantities and denoting their derivatives along $z$ or $z'$
by a subscript $_z$ or $_{z'}$ respectively, we find that the
differentials of the basis 1-forms read as follows.
\beq
&\der E=2\g E\dz F-2\g_{z'}\La\dz E-2\om_{z'}E\dz F'+2\om E\dz K\nonumber\\
&+2\Phi K\dz F+4\Psi K\dz F'\label{eq:dere}\\
&+\Phi_{z'z'}\La\dz F'+[\frac{1}{3}\Psi_{zz}+\frac{1}{12}R+\phi_z]\La\dz F
\nonumber\\
&+[\Psi_z+\phi][K\dz\La-F\dz F']+\Phi_{z'}[K\dz\La+F\dz F'],\nonumber
\eeq
\beq                                                   
&\der F=E\dz K\nonumber\\
&+[\g_{z'}-2\g'_z-\om'_{zz'}-a_{zz'}]F\dz\La+\g'_{z'}F'\dz\La\label{eq:frz}\\
&+[\om'_{z'}+\g']K\dz\La+[\g'+\om_{z'}+a_{z'}]F'\dz F+[\om -\om'+a]F\dz K,
\nonumber
\eeq
\beq
&\der\La=F\dz E'+F'\dz E\nonumber\\
&+[\g-\om'_z-a_z]\La\dz F+[\g'-\om_{z'}-a_{z'}]\La\dz F'+[\om+\om'+a]
\La\dz K\label{eq:larz},
\eeq
\beq
&\der K=\g_{zz'}F\dz\La+\g'_{zz'}F'\dz\La+
[\g_{z'}+\g'_z+a_{zz'}]\La\dz K+[\om'_z-\om_z-2\g-a_z]K\dz F\nonumber\\
&+[\g'_z+a_{zz'}-\g_{z'}+\om'_{zz'}-\om_{zz'}]F\dz F'
+[\om_{z'}-\om'_{z'}-2\g'-a_{z'}]K\dz F'.
\eeq
Differentials of $E'$ and $F'$ are obtained from the above equations
by means of transformations (\ref{eq:trpr}).\\
To calculate differentials of $\bar{E}$ and $\bar{E'}$ we need to know
what are the reality conditions for the null cotetrad (\ref{eq:nt}) on
$\cal M$. If we assume the reality conditions (\ref{eq:reco}) then
we find that
\be
\bar{F}=(1-\eps z\bar{z})F'-\eps\bar{z}z'F-\eps\bar{z}\La+
(\bar{z}-z'+\eps\bar{z}(zz'-\eps))K,
\ee
\beq
&\bar{\La}=(1-|\eps|+\eps\bar{z}\bar{z}')\La\nonumber\\
&+(~(1-|\eps|+\eps\bar{z}\bar{z}')z-\bar{z}'~)F'
+(~(1-|\eps|+\eps\bar{z}\bar{z}')z'-\bar{z}~)F\\
&+[~(z'-\bar{z})(\bar{z}'-z)-\eps (\eps -\bar{z}\bar{z})(\eps-zz')~]K,
\nonumber
\eeq
\be
\bar{K}=(1-|\eps|+\eps zz')K-\eps z'F-\eps zF'-\eps\La.
\ee
The corresponding formula for $\bar{F}'$, may be obtained from
$\bar{F}$ by applying (\ref{eq:trpr}).\\
Now, using the above expressions and (\ref{eq:dere}), (\ref{eq:recocon})
we easily get formulae for $\der\bar{E}$ and $\der\bar{E'}$. These, in
particular, imply that
\beq
&\der\bar{E}\dz\bar{E}\dz\La\dz F'\dz E'\nonumber\\
&=\{~(1-|\eps|)~[~4\bar{\Psi}+2(\bar{\Psi}_{\bar{z}}+\bar{\phi})(\bar{z}-z')
+(\frac{1}{3}\bar{\Psi}_{\bar{z}\bar{z}}+\frac{1}{12}R+\bar{\phi}_{\bar{z}})
(\bar{z}-z')^2~]\label{eq:prcp}\\
&+\eps ~[~-2\bar{\Phi}+2\bar{\Phi}_{\bar{z}'}z'(\eps -z'\bar{z}')-
\bar{\Phi}_{\bar{z}'\bar{z}'}(\eps -z'\bar{z}')^2~]~\}~
K\dz F\dz \bar{E}\dz\La\dz F'\dz\bar{E}'.\nonumber
\eeq
As usual, the primed counterpart of (\ref{eq:prcp}) is obtained by applying
(\ref{eq:trpr}).\\

\section{Appendix D}
Differentials of the basis 1-forms on $\Pen$ are obtained by
using their representation (\ref{eq:nt}), (\ref{eq:fbb}). 
\beq
&\der F= E\dz K+[\G_{212}+z(\G_{132}-\G_{213})-z^2\G_{133}]F\dz P\nonumber\\
&+[\G_{213}-\G_{231}+z\G_{133}]F\dz B+
[\G_{214}+z(\G_{134}-\G_{211})-z^2\G_{131}]F\dz K\nonumber\\
&+[\G_{234}-z\G_{231}]B\dz K+
[\G_{232}-z\G_{233}]B\dz P\nonumber
\eeq
\beq
&\der B= E\dz P+[-\G_{434}+z(\G_{431}-\G_{314})+z^2\G_{311}]B\dz K\nonumber\\
&+[-\G_{431}+\G_{413}-z\G_{311}]B\dz F+
[-\G_{432}+z(\G_{433}-\G_{312})+z^2\G_{313}]B\dz P\nonumber\\
&+[-\G_{412}+z\G_{413}]F\dz P+
[-\G_{414}+z\G_{411}]F\dz K\nonumber
\eeq
\beq
&\der E=[\G_{213}+\G_{343}+2z\G_{133}]E\dz B\nonumber\\
&+[\G_{341}+\G_{211}+2z\G_{131}]E\dz F\nonumber\\
&+[\G_{212}+\G_{342}+z(2\G_{132}-\G_{213}-\G_{343})-2z^2\G_{133}]E\dz P
\nonumber\\
&+[\G_{214}+\G_{344}+z(2\G_{134}-\G_{211}-\G_{341})-2z^2\G_{131}]E\dz K
\nonumber\\
&+[\frac{1}{2}S_{44} -zS_{41} + \frac{1}{2} z^2 S_{11}]K\dz F\nonumber\\
&+[\Psi_0-4\Psi_1z+6\Psi_2z^2
-4\Psi_3z^3+\Psi_4z^4]K\dz P\nonumber\\
&+[\Psi_2+\frac{1}{12}R+\frac{1}{2}\phi_1-
z(2\Psi_3+\frac{1}{2}\phi_0)+z^2\Psi_4]K\dz P\nonumber\\
&+[\frac{1}{2}S_{22}-zS_{32}+\frac{1}{2}z^2S_{33}]P\dz B\nonumber\\
&+[-\Psi_1-\frac{1}{4}\phi_2+z(3\Psi_2+\frac{1}{2}\phi_1)-z^2(3\Psi_3+
\frac{1}{4}\phi_0)-z^3\Psi_4](P\dz F-K\dz B)\nonumber\\
&+[\frac{1}{2}S_{42}-\frac{1}{2}z(S_{12}+S_{34})
+\frac{1}{2}z^2S_{31}](P\dz F+K\dz B)
\nonumber
\eeq
\beq
&\der K= \G_{313}B\dz F+\G_{323}B\dz P+
[\G_{321}-\G_{312}+z\G_{313}]F\dz P\nonumber\\
&+[\G_{341}-\G_{314}]F\dz K+[\G_{343}-z\G_{313}]B\dz K\nonumber\\
&+[-\G_{324}+\G_{342}+z(-\G_{312}+\G_{321}-\G_{343})+z^2\G_{313})]P\dz K
\nonumber
\eeq
\beq
&\der P= [-\G_{132}+\G_{123}]P\dz B+[\G_{121}-z\G_{131}]P\dz F\nonumber\\
&+\G_{131}B\dz F+\G_{141}K\dz F+
[\G_{134}-\G_{143}-z\G_{131}]B\dz K\nonumber\\
&+[\G_{124}-\G_{142}+z(\G_{143}-\G_{134}-\G_{121})+z^2\G_{131}]P\dz K\nonumber
\eeq   
The complex conjugo of $\der E$ is easily calculable using (\ref{eq:reco})
and (\ref{eq:recocon}).\\
Finally, differentials of the basis 1-forms on $\Pen'$ follows from the above
by applying (\ref{eq:trpr}).
\section{Appendix E}
In this Appendix we study the properties of the basis 1-forms on $\pp$ from the
point of view of the Weyl transformations
$$(g, A)\rightarrow (\hat{g}, \hat{A})=
({\rm e}^{2\varphi}g,~ A-2{\rm d}\varphi).~~~~~~~~~~~~~~~~~~~~~~~(*)$$
We start with the local representation (\ref{eq:fonpp}) of the forms. Then we
note that $E$ has the following decomposition onto the Levi-Civita and
$A$ dependent part
\footnote{Analogous formulae may also be obtained for $E'$.}.
$$E=E^{LC}+E^A,$$
where
$$E^{LC}=dz-\gamma^3_{~2} 
+z(\gamma^1_{~1}-\gamma^4_{~4})+z^2\gamma^2_{~3},$$
and
$$E^A=\frac{1}{2}(A_2-zA_3)\La+\frac{1}{2}(zA_1+z'A_2-zz'A_3-A_4)F.$$
We now represent the Weyl transformations as such transformations that
change coefficients of the metric and do not affect the null cotetrad.
Thus we have
$$(g_{ij},A_i)\rightarrow (\hat{g}_{ij},\hat{A}_i)=({\rm e}^{2\varphi}g_{ij},
A_i-2\varphi_{|i}),$$
where the subscript $_{|i}$ means a derivative along the null tetrad vector
$e_i$. Now, it is easy to see that under the above transformations
the forms $E^{LC}$ and $E^A$ transform as follows.
$$\hat{E}^{LC}=E^{LC}+(\varphi_{|2}-z\varphi_{|3})\La +
(z\varphi_{|1}+z'\varphi_{|2}-zz'\varphi_{|3}-\varphi_{|4})F,$$
$$\hat{E}^A=E^A-(\varphi_{|2}-z\varphi_{|3})\La -
(z\varphi_{|1}+z'\varphi_{|2}-zz'\varphi_{|3}-\varphi_{|4})F.$$
This shows that $E$ is invariant under the Weyl transformations. The
same is also true for $F$, $F'$ and $\La$. This, in particular, implies
the invariance of the (anti)-self-duality equations (\ref{eq:weyl}) 
under the Weyl transformations.
The invariance of the Einstein-Weyl equations (\ref{eq:ric})
under these transformations
also follows.\\
Finally we comment on the purely Riemannian case. In this
case we have a given metric $g$ and $A=0$. Since our twistor bundles are
constructed out of null objects then it is reasonable to ask how the
constructions change under the conformal transformations
$g\rightarrow\hat{g}={\rm e}^{2\varphi}g$ of the metric. It follows from
the above transformations of $E^{LC}$
and from the conformal invariance of $\La$, $F$ and $F'$
that although the twistor
bundles for any metric from a given conformal class are the same, the
horizontal spaces for different base metrics are different. This difference
is not essential for the system of (anti)-self-duality equations
(\ref{eq:integr1}), (\ref{eq:integr2}),
(\ref{eq:weyl}) which
is invariant under the conformal transformations, 
but it is essential for the Einstein equations (\ref{eq:ric}). In this
latter case, we need to pick up a particular metric on $\cal M$ and then
use it to define the forms $E$, $\La$, $F'$. Using them we can encode
the Einstein equations for $g$ in $\pp$.

\end{document}